\documentclass[a4paper,12pt]{article}

\newenvironment{proof}{\noindent {\bf Proof:}}{\hfill $\Box$}

\newtheorem{theorem}{Theorem}
\newtheorem{lemma}{Lemma}
\newtheorem{proposition}{Proposition}
\newtheorem{corollary}{Corollary}
\newtheorem{algorithm}{Algorithm}

\newtheorem{definition}{Definition}

\newtheorem{remark}{Remark}

\usepackage{booktabs}
\usepackage{color}

\usepackage{algorithm}
\usepackage[noend]{algpseudocode}
\makeatletter
\def\BState{\State\hskip-\ALG@thistlm}
\makeatother
\usepackage{algorithmicx}

\usepackage{mathrsfs}

\usepackage{amsmath}
\usepackage{amssymb}
\usepackage{amsfonts}
\usepackage{graphicx}
\usepackage{hyperref}

\usepackage{nicefrac}

\usepackage[normalem]{ulem}
\usepackage[dvipsnames]{xcolor}

\newcommand{\R}{\mathbb{R}}
\newcommand{\C}{\mathcal{C}}
\newcommand{\N}{\mathbb{N}}

\begin{document}

\title{\bf Sparsity structures for Koopman operators}

\maketitle

\begin{center}
\author{Corbinian Schlosser$^{1}$},
\author{Milan Korda$^{\text{1},\text{2}}$}
\end{center}

\begin{center}
\date{\small \today}
\end{center}

\footnotetext[1]{CNRS; LAAS; 7 avenue du colonel Roche, F-31400 Toulouse; France. {\tt cschlosser@laas.fr, korda@laas.fr}}
\footnotetext[2]{Faculty of Electrical Engineering, Czech Technical University in Prague,
Technick\'a 2, CZ-16626 Prague, Czech Republic.}

\begin{abstract}
  	We present a decomposition of the Koopman operator based on the sparse structure of the underlying dynamical system, allowing one to consider the system as a family of subsystems interconnected by a graph. Using the intrinsic properties of the Koopman operator, we show that eigenfunctions for the subsystems induce eigenfunctions for the whole system. The use of principal eigenfunctions allows to reverse this result. Similarly for the adjoint operator, the Perron-Frobenius operator, invariant measures for the dynamical system induce invariant measures of the subsystems, while constructing invariant measures from invariant measures of the subsystems is less straightforward. We address this question and show that under necessary compatibility assumptions such an invariant measure exists. Based on these results we demonstrate that the a-priori knowledge of a decomposition of a dynamical system allows for a reduction of the computational cost for data driven approaches on the example of the dynamic mode decomposition.
\end{abstract}

\textbf{keywords:} Koopman theory, dynamical system, model reduction, sparsity, invariant measure, dynamic mode decomposition, sum-of-squares

\section{Introduction}
Given a dynamical system $x^+ = f(x)$ for a function $f:X \rightarrow X$, the Koopman operator for this dynamical system is defined by $Tg := g \circ f$ for functions $g:X \rightarrow \mathbb{C}$ in a suitable function space. Importantly, $T$ is a linear operator. Exploiting this linear structure has many applications; see for example in~\cite{AppliedKoopmanism},~\cite{MezicErgodicPartition},~\cite{OpTheoErgodicTheory},~\cite{Viktoria} or~\cite{MilanMPC} to name just a few of them.
	
We address exploitation of sparse structures following the concepts of sparsity from~\cite{SparsePaper}, which is strongly motivated by~\cite{Chen}. We use the same concept of sparsity, i.e., the property that the system decouples into subsystems interconnected according the underlying so-called sparsity graph, and describe structures for the dynamical system that inherit such sparsity. For the Koopman operator this means that we aim to describe objects such as eigenfunctions by corresponding objects for smaller subsystems induced by the sparsity patterns considered.
	
In terms of exploiting structures our work is related to~\cite{Symmetry} where it is shown that symmetry translates to operators commuting with the Koopman operators. In this text we show that the considered sparse structures allow to \emph{intertwine} the Koopman operator $T$ of the whole system with the Koopman operators $T^I$ for corresponding smaller systems, i.e. there exists an injective linear operator $C$ such that
\begin{equation*}
    CT = T^IC.
\end{equation*}
The operator $C$ appears naturally and is the composition operator for the projection onto the subsystem.

Both the approaches, via symmetry \cite{Symmetry} and sparsity in this text, formulate properties of the dynamical system as properties of the Koopman operator. This follows the classical idea of Koopman theory (see for example~\cite{KariKuester}), namely to translate properties of $f$ or its corresponding dynamical system into functional analytic properties of the Koopman operator -- and vice versa. In many cases this procedure allows to apply spectral theory. This does not just induce a better understanding of the evolution of the dynamical system but also gives strong theoretical foundations for many results, as it lies at the core of several applications. For instance for (extended) dynamic mode decomposition ((E)DMD) the importance of the spectral theory for the Koopman operator was demonstrated in~\cite{MilanConvergence} and~\cite{MilanDMD} where convergence of the DMD was proved using results on spectral projections. Other examples, where spectral theory for the Koopman operator induces results for the dynamical system, can be found in~\cite{MezicErgodicPartition} or in~\cite{MezicMauroy} among others.
	
The price one has to pay for obtaining a linear representation of (nonlinear) dynamical systems is that the state space for the Koopman operator is infinite dimensional. This comes along with high computational costs when trying to extract the information contained in the Koopman operator on the function space. Hence reduction techniques are important tools for computational methods. While for example in~\cite{KernelMethodsDMD} a reduction of complexity by using kernel methods is suggested, we follow the approach of decomposing the Koopman operator into Koopman operators of smaller systems. We demonstrate this by numerical examples of a reduced DMD and sparse computation of invariant measures. 

\section{Notations}
The non-negative real numbers are denoted by $\R_+$ and the natural numbers are denoted by $\N$. The cardinality of a set $A$ is denoted by $|A|$. A set $I$ always denotes a subset of $\N$ and if $I \subset \{1,\ldots,n\}$ the corresponding projections in $\R^n$ onto the canonical coordinates indexed by $I$ is denoted by $\Pi_I:\R^n \rightarrow \R^{|I|}$, i.e. $\Pi_I(x_1,\ldots,x_n) = (x_i)_{i \in I}$. For a matrix $A \in \R^{n \times m}$ we denote by $A^T$ its transpose. For a function $f:\R^n \rightarrow \R^n$ we denote by $f_I := \Pi_I \circ f$ the restriction of $f$ to its coordinates indexed by $I$. We denote by $\|x\|$ the euclidean norm of an element $x \in \R^n$. The space of polynomials on $\R^n$ is denoted by $\R[x]$ or $\R[x_1,\ldots,x_n]$ and the space of polynomials of degree at most $d \in \N$ by $\R[x]_d$. For a polynomial $g$ we denote by $\mathrm{deg}(g)$ its degree. For a subset $A$ of a topological space $X$ we denote by $\overline{A}$ its closure. The space of continuous real-valued functions on a topological space $X$ is denoted by $\C(X)$, and when $X$ is compact we equip implicitly $C(X)$ with the supremum norm. For continuously differentiable functions $g$ we denote their derivative by $Dg$ and the second derivative by $D^2g$. The space of signed Borel measures, is denoted by $\mathcal{M}(X)$ and we interpret $\mathcal{M}(X)$ as the dual of $\C(X)$ when $X$ is compact. The cone of non-negative measures on $X$ is denoted by $\mathcal{M}(X)_+$. The action of a measure $\mu \in \mathcal{M}(X)$ on a continuous function $g \in \C(X)$ by integration is denoted by $\langle g,\mu \rangle := \int\limits_X g \; d\mu$. For two topological spaces $X,Y$ and a Borel measureable map $h:X \rightarrow Y$ and a measure $\mu \in \mathcal{M}(X)$ let $h_\#\mu$ be the push forward measure defined by $(h_\# \mu) (A):= \mu(h^{-1}(A))$ for all Borel sets $A \subset Y$. For a linear operator $V$ we denote its dual by $V^*$.

\section{Setting and preliminary definitions}

In this article we consider dynamical systems on $\R^n$ induced by an ordinary differential equation \begin{equation}\label{EqnODE}
	\dot{x} = f(x), \;\; x(0) = x_0 \in \R^n.
\end{equation}
for a Lipschitz continuous vector field $f:\R^n \rightarrow \R^n$. The corresponding semiflow/solution map is denoted by $\varphi_t$ for $t \in \R_+$.

\begin{remark}\label{RemarkDiscrete}
    We focus on continuous time systems but the treatment of discrete time dynamical systems is the same -- we only have to replace the continuous time objects by their discrete time analogues.
\end{remark}

We focus on the notion of a subsystem and corresponding properties of the Koopman and Perron-Frobenius operators.

For a subset of indices $I \subset \{1,\ldots,n\}$ we denote the corresponding projections in $\R^n$ along the canonical coordinates indexed by $I$ by $\Pi_I:\R^n \rightarrow \R^{|I|}$, i.e. $\Pi_I(x_1,\ldots,x_n) = (x_i)_{i \in I}$

\begin{definition}\label{def:Subsystem}
    For $I \subset \{1,\ldots,n\}$ we call $(I,f_I)$ a subsystem or a subsystem induced by $I$ if $f_I := \Pi_I \circ f$  only depends on the states index by $I$.
\end{definition}

\begin{remark}[Coordinate-free formulation]
The definition of a subsystem is coordinate dependent. A coordinate free formulation can be developed within the framework of the so-called factor systems. This abstract treatment is currently not amenable for computation and is briefly described in Section \ref{subsec:CoordinateFreeFormulation}.
\end{remark}

The idea of a subsystem $(I,f_I)$ is that we can treat it as a lower dimensional dynamical system, namely on $\R^{|I|}$ instead of $\R^n$. To see this we view $f_I$ as a (Lipschitz) vector field on $\R^{|I|}$. To do so we identify $f_I$ with the map from $\R^{|I|}$ to $\R^{|I|}$ given by $(x_i)_{i \in I} \mapsto (f_i(x))_{i \in I}$ where $x = (x_1,\ldots,x_n) \in \R^n$ is any vector with $\Pi_I(x) = (x_i)_{i \in I}$, for example $x_j = 0$ whenever $j \notin I$. Due to the condition in the Definition of a subsystem this map is indeed independent of the specific choice of such $x$. We then have
\begin{equation}\label{equ:def:subsystem}
    f_I \circ \Pi_I = \Pi_I \circ f.
\end{equation}
The semiflow induced by $f_I:\R^{|I|} \rightarrow \R^{|I|}$ is denoted $\varphi_t^{I}:\R^{|I|} \rightarrow \R^{|I|}$ for $t \in \R_+$ and by virtue of (\ref{equ:def:subsystem}) it satisfies 
\begin{equation}\label{equ:FlowCommuteSubsystem}
    \varphi_t^I \circ \Pi_I = \Pi_I \circ \varphi_t.
\end{equation}

The easiest non-trivial example of a dynamical system inducing subsystems is taken from~\cite{Chen} and has the following form
\begin{align}\label{equ:Sparsex0x1x2}
	\dot{x}^1 & = & f_1(x^1)\notag\\
	\dot{x}^2 & = & f_2(x^1,x^2)\\
	\dot{x}^3 & = & f_3(x^1,x^3)\notag
\end{align}
for $x^1 \in \R^{n_1},x^2 \in \R^{n_2},x^3 \in \R^{n_3}$ for $n_1,n_2,n_3 \in \N$ and maps $f_1:\R^{n_1}\rightarrow \R^{n_1}, f_2: \R^{n_1 + n_2} \rightarrow \R^{n_2},f_3: \R^{n_1 + n_3} \rightarrow \R^{n_3}$. In this case, the state space is $\R^{n_1 + n_2 + n_3}$ and $I_1 = \{1,\ldots,n_1\}$, $I_2 = \{1,\ldots,n_1 + n_2\}$ and $I_3 = \{1,\ldots,n_1,n_1 + n_2+1,\ldots,n_1+ n_2 + n_3\}$ induce subsystems. That means in our setting we have $f:\R^{n_1 + n_2 + n_3} \rightarrow \R^{n_1 + n_2 + n_3}$ with $f(x) = (f_1(\Pi_{I_1}(x)),f_2(\Pi_{I_2}(x)),f_3(\Pi_{I_3}(x)))$.

\paragraph{Sparsity graph} An essential tool to describe the sparsity of $f$ is its \emph{sparsity graph} (see \cite{SparsePaper} for a treatment of sparsity graphs and subsystems).

\begin{definition}[Sparsity graph]\label{def:sparsitygraph}
    Let $(x_1,\ldots,x_n)$ denote a variable in $\R^n$. The sparsity graph associated to $f = (f_1,\ldots,f_n):\R^n \rightarrow \R^n$ is defined by:
    \begin{enumerate}
        \item The set of nodes is $(x_1,\ldots,x_n)$.
        \item For $i \neq j$ there is an edge $(x_i,x_j)$ if the function $f_j$ depends on the variable $x_i$.
    \end{enumerate}
\end{definition}

Subsystems $(I,f_I)$ can be described, in terms of the sparsity graph, by exactly those collections $(x_i)_{i \in I}$ of nodes such that all incoming edges to any $x_i$ for $i \in I$ starts at a node $x_j$ for some $j \in I$. Finding subsystems can then be done by simple graph methods applied to the sparsity graph of $f$ (\cite{SparsePaper}).

\begin{figure}[!h]
\begin{center}
\includegraphics[width=130mm]{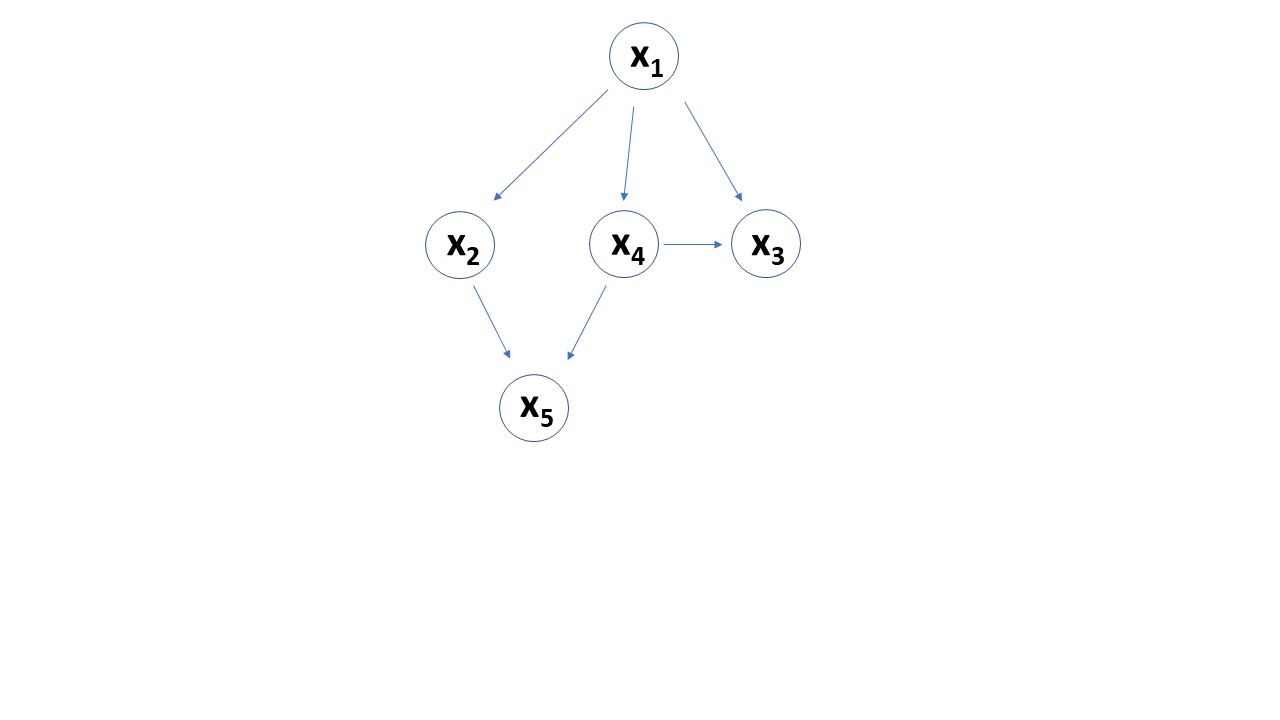}
\end{center}
\vspace{-30mm}
\caption{\footnotesize{Sparsity graph for the states $x_1,\ldots,x_5$ and dynamics given by $\dot{x}_1 = f_1(x_1)$, $\dot{x}_2 = f_2(x_1,x_2)$, $\dot{x}_3 = f_3(x_1,x_3,x_4)$, $\dot{x}_4 = f_4(x_1,x_4)$ and $\dot{x}_5 = f_5(x_2,x_4)$ for vector fields $f_1,\ldots,f_5$, i.e. $f(x) = (f_1(x_1),f_2(x_1,x_2),f_3(x_1,x_3,x_4),f_4(x_1,x_4),f_5(x_2,x_4))$. The subsystems are induced by $I_1 = \{1\}, I_2 = \{1,2\}, I_3 = \{1,4\}, I_4 = \{1,3,4\}, I_5 = \{1,2,4,5\}, I_6 = \{1,2,4\}, I_7 = \{1,2,3,4\}$ and $I_8 = \{1,2,3,4,5\}$.}}
\label{fig:Example}
\end{figure}

Figure \ref{fig:Example} displays the sparsity graph and its relation to subsystems for a less trivial example.

In this work, we are interested in the case when the dynamics are constrained to a set $X \subset \R^n$. For subsystem $(I,f_I)$ we always consider the corresponding constraint set $\Pi_I(X)$.

We say $X \subset \R^n$ is positively invariant if for all $x \in X$ the trajectory $\varphi_t(x)$ exists for all $t \in \R_+$ and $\varphi_t(x) \in X$ for all $t \in \R_+$.

\textbf{Assumption:} $X$ is positively invariant. \label{assumption:XPosInv}

Note that if $X$ is positively invariant and $I$ induces a subsystem then $\Pi_{I}(X)$ is positively invariant.

We need a last notion concerning sparsity which plays an important role in Theorem \ref{thm:GluedInvariantMeasure}. Therefore, we will assume that not just the dynamics allow a projection onto a subset of all coordinates but also the constraint set $X$ has a factoring property that is compatible with the concept of subsystems.

\begin{definition}\label{def:Factorization}
	We say $X \subset \R^n$ factors with respect to a family of subsystems induced by $I_1,\ldots,I_N$ if $\bigcup\limits_{k = 1}^N = \{1,\ldots,n\}$ and
	\begin{equation*}
		X = \{x \in \R^n: \Pi_{I_k} (x) \in \Pi_{I_k}(X) \text{ for } k = 1,\ldots,N\}.
	\end{equation*}
\end{definition}

For subsystems induced by $I_1 = \{1,\ldots,n_1\}, I_2 = \{n_1 + 1,\ldots,n_2\}, \ldots, I_N = \{n_{N-1} + 1,\ldots,n\}$ this means that $X$ is of the form
\begin{equation}\label{equ:XProduct}
X = X_1 \times \ldots \times X_N
\end{equation}
with $X_j \subset \R^{|I_j|}$ for $j = 1,\ldots,N$. In case for the example from (\ref{equ:Sparsex0x1x2}), that means $X$ factors with respect to the subsystems $I_1 = \{1,2\}$ and $I_2 = \{1,3\}$ if it is of the form $X = X_1 \times X_2 \times X_3$ for $X_i \subset \R^{n_i}$ for $i = 1,2,3$. We refer to \cite{SparsePaper} for more details on the notion of factoring with respect to subsystems.

\paragraph{Koopman and Perron-Frobenius operators}

Now we define Koopman and Perron-Frobenius operators corresponding to a dynamical system.

\begin{definition}\label{Def:Koopman}
	Let $X$ be positively invariant for (\ref{EqnODE}) and $t \in \R_+$. The Koopman operator $T_t:\C(X)\rightarrow \C(X)$ is given by
	\begin{equation}\label{equ:DefKoopman}
		T_tg := g \circ \varphi_t.
	\end{equation}
	The Perron-Frobenius operator $P_t:\mathcal{M}(X) \rightarrow \mathcal{M}(X)$ is given by $P_t = (\varphi_t)_\#$, i.e.
	\begin{equation}\label{equ:DefPerronFrobenius}
		P_t\mu(A) = \mu(\varphi_t^{-1}(A))	
	\end{equation}
	for any Borel set $A \subset X$.
\end{definition}

It is a standard result that if $X$ is compact the Perron-Frobenius operator is the adjoint of the Koopman operator, i.e.
\begin{equation}\label{Eq:AdjointCondition}
	\int\limits_X g \circ \varphi_t \; d\mu = \int\limits_X g \; dP_t\mu,
\end{equation}
and the family of operators $(T_t)_{t \in \R_+}$ is a strongly continuous semigroup of contractive linear algebra homomorphisms, i.e. $T_t$ is a linear operator for all $t \in \R_+$, $T_{t+s} = T_tT_s$ for all $t,s \in \R_+$, for all $g \in \C(X)$ we have $T_tg \rightarrow g$ uniformly as $t \rightarrow 0$, $\max\limits_{x \in X} |T_tg(x)| \leq \max\limits_{x \in X} |g(x)|$ and $T_t(g_1 \cdot g_2) = T_tg_1 \cdot T_t g_2$ for all $t \in \R_+$ and $g_1,g_2 \in \C(X)$.

\section{Sparse properties of the Koopman operator induced by subsystems}

%

In the following we will be mostly interested in eigenfunctions and eigenmeasures, i.e. eigenvectors of the Koopman and Perron-Frobenius operator respectively. We say $g \in \C(X)$ respectively $\mu \in \mathcal{M}(X)$ is an eigenfunction respectively eigenmeasure with eigenvalue $\lambda \in \mathbb{C}$ of the Koopman respectively Perron-Frobenius operator if $g \neq 0$ respectively $\mu \neq 0$ and for all $t \in [0,\infty)$ we have $T_tg = e^{\lambda t}g$ and $P_t \mu = e^{\lambda t} \mu$ respectively. Of particular interest for the Perron-Frobenius operator are invariant measures, i.e. non-negative measures $\mu \in \mathcal{M}(X)_+$ such that $P_t\mu = \mu$.\\
Both, eigenfunctions of the Koopman operator and eigenmeasures of the Perron-Frobenius operator, have important applications. For example they allow to generalize the concept of ``diagonalizing" the dynamics or give insight into ergodic properties of the dynamical system (see for example \cite{AppliedKoopmanism} and \cite{OpTheoErgodicTheory}).

Before stating the following proposition we want to fix some notation. The Koopman operator for the subsystem induced by $I$ with corresponding constraint set $\Pi_I(X)$ is denoted by $T_t^I$. Hence $T_t^I$ acts on $\C(\Pi_I(X))$; the Perron-Frobenius operator for the subsystem is denoted by $P_t^I$. We say an operator $V:W\rightarrow Z$ intertwines an operator $S:Z \rightarrow Z$ with $T:W \rightarrow W$ if
\begin{equation}\label{equ:IntertwineProjectionKoopman}
	S V = V T.
\end{equation}

The next Proposition is a direct consequence from the fact that the operators $V_I:\C(\Pi_I(X)) \rightarrow \C(X)$ and $V_I^*:M(X) \rightarrow \mathcal{M}(\Pi_I(X))$ given by
\begin{equation}\label{eq:IntertwineOp}
V_Ig := g \circ \Pi_I \text{ and } V_I^* \mu := (\Pi_I)_\#\mu
\end{equation}
intertwine $T_t$ and $T_t^I$ respectively $P_t^I$ and $P_t$ for all $t \in \R_+$ (\cite{OpTheoErgodicTheory} page 233). Even though the statement of Proposition \ref{propEigenvector} is known
and follows immediately from (\ref{eq:IntertwineOp}) we state its proof because it illustrates how we apply the relation (\ref{equ:FlowCommuteSubsystem}) to the Koopman and Perron-Frobenius operator.

\begin{proposition}\label{propEigenvector}
	Let $X$ be positively invariant and $I$ induce a subsystem. Then
	\begin{enumerate}
	    \item $V_I$ is injective and intertwines $T_t$ and $T_t^I$ and $V^*_I$ intertwines $P_t^I$ and $P_t$ for all $t \in \R_+$. If $X$ is compact then $V_I^*$ is the dual of $V_I$ and surjective.
		\item If $g \in \C(\Pi_I(X))$ is an eigenfunction with eigenvalue $\lambda$ for the Koopman operator $T_t^I$ for the corresponding subsystem then $\hat{g} \in \C(X)$ defined by $\hat{g}:= g \circ \Pi_I$ is an eigenfunction  with eigenvalue $\lambda$ of the Koopman operator $T_t$ for the whole system.
		\item If $\mu \in \mathcal{M}(X)$ is an eigenmeasure with eigenvalue $\lambda$ of the Perron-Frobenius operator $P_t$, so is the push forward measure of $\mu$ by $\Pi_I$, i.e. $\mu_I := (\Pi_I)_\# \mu$, an eigenmeasure with eigenvalue $\lambda$ for the Perron-Frobenius operator for the subsystem $P_t^I$.
	\end{enumerate}
\end{proposition}

\begin{proof}
    Let $g \in \C(\Pi_I(X))$, $\mu \in \mathcal{M}(X)$ and $t \in \R_+$.
	\begin{enumerate}
		\item Because $\Pi_I:X \rightarrow \Pi_I(X)$ is clearly surjective it follows that $V_I$ is injective. To check this let $g \in \C(\Pi_I(X))$ with $0 = V_Ig = g\circ \Pi_I$, i.e. $0 = g(\Pi_I(x))$ for all $x \in X$, hence $g = 0$ and $V_I$ is injective. From (\ref{equ:FlowCommuteSubsystem}) we get
		\begin{eqnarray}\label{eq:IntertwineCalculation}
			T_t V_I g & = & T_t (g\circ \Pi_I) = (g \circ \Pi_I) \circ \varphi_t = g \circ (\Pi_I \circ \varphi_t)\notag\\
								& = & g \circ (\varphi^I_t \circ \Pi_I) = (g \circ \varphi^I_t) \circ \Pi_I = V_I (g \circ \varphi_t^I) = V_I T_t^Ig.
		\end{eqnarray}
		When $X$ is compact $V_I^*$ is the taking the adjoint of $V_I$ and the adjoint of $T_tV_I = V_I T_t^I$ gives \begin{align*}
		   V_I^*P_t = V_I^*T_t^* = (T_tV_I)^* = (V_IT_t^I)^* = P_t^I V_I^*.
		\end{align*}
		For the case when $X$ is not compact we can also argue in an analogue way to (\ref{eq:IntertwineCalculation}) to obtain that $V_I^*$ intertwines $P_t^I$ and $P_t$. For the proof that $V_I$ is surjective when $X$ is compact we refer to \cite{OpTheoErgodicTheory} page 208.
		\item Assume $g$ is an eigenfunction for $T_t^I$. Then for all $t \in \R_+$ we have $e^{\lambda t} g = T_t^I g = g \circ \varphi_t^I$ and by (\ref{equ:FlowCommuteSubsystem}) we get for $\hat{g}$
		\begin{eqnarray*}
			T_t \hat{g} & = & \hat{g} \circ \varphi_t = (g \circ \Pi_I) \circ \varphi_t = g \circ (\Pi_I \circ \varphi_t)\\
							 & = & g \circ (\varphi_t^I \circ \Pi_I) = (g \circ \varphi_t^I) \circ \Pi_I = e^{\lambda t} g \circ \Pi_I = e^{\lambda t} \hat{g}.
		\end{eqnarray*}
		\item Similarly for the dual operator we have for eigenmeasures $\mu$ of $P_t$ with eigenvalue $e^{\lambda t}$
		\begin{align*}
			P_t^I \mu_I = & (\varphi_t^I)_\# (\Pi_I)_\# \mu = (\varphi_t^I \circ \Pi_I)_\#\mu = (\Pi_I \circ \varphi_t)_\mu\\
						= & (\Pi_I)_\# (\varphi_t)_\# \mu = e^{\lambda t} (\Pi_I)_\# \mu = e^{\lambda t} \mu_I.
		\end{align*}
	\end{enumerate}
\end{proof}

The converse question of constructing eigenfunctions for the subsystem from eigenfunctions for the whole system and analogue constructing eigenmeasures for the whole system from eigenmeasures for the subsystem is less straightforward. 
We treat that problem in Sections \ref{Subsec:InvMeaure} and \ref{Subsec:Eigenfunction} in Theorems \ref{thm:GluedInvariantMeasure} and \ref{thm:PrincipalEigenvalues}.

\subsection{Construction of eigenmeasures from eigenmeasures of subsystems}\label{Subsec:InvMeaure}

In this section we present that for given invariant measures for subsystems satisfying necessary compatibility conditions we can construct (or glue together) those measures to obtain an invariant measure for the whole system.

Even if we will not give an explicit construction of how to ``glue together" the invariant measures for the subsystems to obtain an invariant measure for the whole system the next theorem states that the invariant measure for the whole system ``arise" from a decomposition into subsystems and their corresponding invariant measures.

\begin{theorem}\label{thm:GluedInvariantMeasure}
	Let $X$ be compact and factor with respect to subsystems $(I_1,f_{I_1}),\ldots, (I_N,f_{I_N})$ according to Definition \ref{def:Factorization}. For $k = 1,\ldots,N$ let an invariant probability measure for the subsystem induced by $I_k$ be given by $\mu_k \in \mathcal{M}(\Pi_{I_k}(X))$. Then there exists an invariant probability measure $\mu \in \mathcal{M}(X)$ such that
	\begin{equation}\label{eq:Marginals}
		(\Pi_{I_k})_\# \mu = \mu_k \text{ for all } k = 1,\ldots,N
	\end{equation}
	if and only if for all $k,l \in \{1,\ldots,N\}$
	\begin{equation}\label{equ:PropNecessaryInvMeasureProjection}(\Pi_{I_k \cap I_l})_\# \mu_k = (\Pi_{I_k \cap I_l})_\#\mu_l
	\end{equation}
\end{theorem}

\begin{proof} Necessity of (\ref{equ:PropNecessaryInvMeasureProjection}) follows from $\Pi_{I_k} \circ \Pi_{I_l} = \Pi_{I_l} \circ \Pi_{I_k} = \Pi_{I_k \cap I_l}$ because for any $k,l \in \{1,\ldots,N\}$ we get
\begin{eqnarray*}
    (\Pi_{I_k \cap I_l})_\# \mu_k & = & (\Pi_{I_l} \circ \Pi_{I_k})_\# \mu_k = (\Pi_{I_l})_\# (\Pi_{I_k})_\# \mu_k = (\Pi_{I_l})_\# \mu_k\\
    & = & (\Pi_{I_l})_\# (\Pi_{I_k})_\# \mu = (\Pi_{I_l} \circ \Pi_{I_k})_\# \mu.
\end{eqnarray*}
Replacing the roles of $k$ and $l$ and using that $\Pi_{I_k}$ and $\Pi_{I_l}$ commute we see that $(\Pi_{I_k \cap I_l})_\# \mu_k = (\Pi_{I_l} \circ \Pi_{I_k})_\# \mu = (\Pi_{I_k \cap I_l})_\# \mu_l$, i.e. (\ref{eq:Marginals}). For the sufficiency part we consider the set
\begin{equation}\label{eq:SetK}
	K := \{\mu \in \mathcal{M}(X)_+: \mu(X) = 1, (\Pi_{I_k})_\# \mu = \mu_k, k = 1,\ldots,N\}
\end{equation}
and will show that it is non-empty, convex, compact (with respect to the weak-$^*$ topology) and $P_t$ invariant for all $t \in \R_+$. The result then follows from the Markov-Kakutani theorem (\cite{OpTheoErgodicTheory}). To show that $K$ is non-empty we recall that by~\cite{AmbrosioOptimalTransport} Lemma 2.1 it is possible to glue two probability measures with coinciding common marginals together ``along the marginal". That means for probability measures $\mu \in \mathcal{M}(X \times Y)$ and $\nu \in \mathcal{M}(X \times Z)$ with $(\Pi_X)_\# \mu = (\Pi_X)_\# \nu$ there exists a probability measure $\gamma \in \mathcal{M}(X\times Y \times Z)$ with $(\Pi_{X \times Y})_\# \gamma = \mu$ and $(\Pi_{X \times Z})_\#\gamma = \nu$. The compatibility condition (\ref{equ:PropNecessaryInvMeasureProjection}) guarantees that the common marginals of the measures $\mu_k$ coincide and we can apply~\cite{AmbrosioOptimalTransport} Lemma 2.1 (inductively) to glue together the measures $\mu_k$. The condition that $X$ factors according to Definition \ref{def:Factorization} assures that the support of such a glued measure is contained in $X$. Hence the set $K$ is non-empty. To check convexity and weak-$^*$ closedness note that the bounded composition operators $V_{I_k}:\C(\Pi_{I_k}(X)) \rightarrow \C(X)$ given by
\begin{equation*}
	V_{I_k} f = f \circ \Pi_{I_k},
\end{equation*}
are linear and bounded, so are their adjoints $V_{I_k}^* = (\Pi_{I_k})_\#:\mathcal{M}(X) \rightarrow \mathcal{M}(\Pi_{I_k}(X))$. From linearity of $V_{I_k}^I$ it follows that $K$ is convex. Boundedness of $V_{I_k}^*$ implies that $K$ is also weak-* closed because bounded operators are continuous with respect to the weak-* topology. Hence the set $W = \bigcap\limits_{k = 1}^N (V_{I_k}^*)^{-1}(\{\mu_k\})$ is weak-* closed and convex as the intersection of weak-* closed convex sets. The constraint $\mu(X) = 1$ implies that $K$ is a (closed, convex) subset of the set of probability measures - hence it is compact with respect to the weak-$^*$ topology. To check that $K$ is $P_t$ invariant for all $t \in \R_+$ let $\mu \in K$, $t \in \R_+$ and $1 \leq k \leq N$. Then $P_t\mu(X) = \mu(\varphi_t^{-1}(X)) = \mu(X) = 1$ and
\begin{eqnarray*}
    (\Pi_{I_k})_\#P_t\mu & = & (\Pi_{I_k})_\# (\varphi_t)_\# \mu = (\Pi_{I_k} \circ \varphi_t)_\# \mu = (\varphi_t^I \circ \Pi_{I_k})_\#\mu\\
    & = & (\varphi_t^{I_k})_\# (\Pi_{I_k})_\# \mu \overset{(\ref{eq:SetK})}{=} P_t^I \mu_k = \mu_k
\end{eqnarray*}
where the last equality holds because $\mu$ is an invariant measure for the subsystem, i.e. $P_t^I\mu_k = \mu_k$. That shows invariance of $K$ with respect to $P_t$ for all $t \in \R_+$. Further, for all $t,s\in \R_+$ the operators $P_t$ and $P_s$ commute, due to
\begin{align}\label{eq:CommuteSemigroup}
P_tP_s = P_{t+s} = P_sP_t.
\end{align}
Because the operators $P_t$ are bounded for all $t \in \R_+$ they are also continuous with respect to the weak-$^*$ topology and we can apply the Markov-Kakutani theorem (\cite{OpTheoErgodicTheory} Theorem 10.1) to the family of operators $(P_t)_{t \in \R_+}$ and the set $K$ from (\ref{eq:SetK}). This gives a measure $\mu \in K$ that satisfies $P_t \mu = \mu$ for all $t \in \R_+$, i.e. $\mu$ is an invariant probability measure with the given marginals $(\Pi_{I_k})_\# \mu = \mu_k$ for $k =1,\ldots,N$.
\end{proof}

\begin{remark}\label{rem:AveragingInvariantMeasure}
    The proof of the Theorem \ref{thm:GluedInvariantMeasure} can be made more explicit by using the technique from the proof of the Markov-Kakutani theorem directly. Namely, in order to find an invariant measure in $K$, we can take any measure $\mu \in K$ and consider the time averages for $T>0$
    \begin{equation*}
	    \mu_T:= \frac{1}{T}\int\limits_0^T P_s \mu \; ds.
    \end{equation*}
    By compactness and convexity of $K$ we have $\mu_T \in K$ for all $T >0$ and there exists a weak-$^*$ converging subsequence $\mu_{T_k}$ for $T_k \rightarrow \infty$ with limit $\mu_\infty \in K$. Then $\mu_\infty$ is an invariant measure.
\end{remark}

Theorem \ref{thm:GluedInvariantMeasure} can be seen as a generalization of a similar result for totally decoupled systems in \cite{OpTheoErgodicTheory} page 208.

In the spirit of averaging the construction of eigenmeasures with eigenvalues $\lambda \neq 0$ can be approached via, for instance, Laplace averages. This is more subtle because the typical problem of existence or boundedness of Laplace averages also appears here. We refer to~\cite{LaplaceAverages} for a treatment of Laplace averages (for eigenfunctions).

\subsection{Eigenfunctions of the Koopman opertor based on eigenfunctions of subsystems}\label{Subsec:Eigenfunction}

For two topological spaces $X$ and $Z$ we have a canonical way of projecting a measure in $\mathcal{M}(X \times Z)$ to measures in $\mathcal{M}(X)$ and $\mathcal{M}(Z)$. For functions $g \in \C(X \times Z)$ it is not so clear how to project $g$ onto $\C(X)$ and $\C(Z)$. The evaluation map $g(\cdot,\cdot) \mapsto g(x_0,\cdot)$ for some $x_0 \in X$ does not send eigenfunctions to eigenfunctions in general. Fortunately, we will see that the so-called principal eigenvalues have a certain decomposition property. The decomposition of principal eigenfunctions will be based on their uniqueness; we use such uniqueness results from \cite{Matthew}.

\begin{definition} For systems with globally exponentially attractive fixed point $x^*$ we call an eigenfunction $g \in \C^1(X)$ principal eigenfunction for the Koopman operator if $Dg(x^*) \neq 0$.
\end{definition} 

The set of eigenfunctions of the Koopman operator can be large. The method from \cite{MilanEigenfunction} provides a possibility of constructing arbitrarily many eigenfunctions. Therefore, principal eigenfunctions are motivated by the important attempt to single out some very characteristic eigenfunctions. The underlying idea is to find a ``basis" of eigenfunctions, in the sense that all other eigenfunctions can be constructed by products and sums of the functions in the ``basis". If $g$ is an eigenfunction with eigenvalue $\lambda \neq 0$ then also $g^r$ is an eigenfunction with eigenvalue $r\lambda$, but if $g$ is differentiable and $r > 1$ then $Dg^r(x^*) = rg^{r-1}(x^*) Dg(x^*) = 0$ because $g(x^*) = 0$ for non constant $g$. So the condition $\nabla g(x^*) \neq 0$ restricts to eigenfunctions that are not obtained by powers of other eigenfunctions, and hence seem to be good candidates for ``basic" eigenfunctions. In the case of linear dynamics, $\dot{x} = Ax$ for a matrix $A$, the principal eigenfunctions are the linear forms $\langle \cdot,w\rangle$ where $w$ is an eigenvectors of $A^T$. This is emphasized by the fact that $\nabla g(x^*)$ is an eigenvector of $Df(x^*)^T$ for any principal eigenfunction $g$ and the corresponding eigenvalue $\lambda$ for $\nabla g(x^*)$ coincides with the eigenvalue for the eigenfunction $g$. This motivates the following definition.

\begin{definition}\label{def:PrincipalEigenfunction}
	We call a function $g \in \C^1(X)$ a principal eigenfunction of the Koopman operator tangential to a vector $ 0 \neq v \in \R^n$ if $\nabla g(x^*) = v$.
\end{definition}

Principal eigenfunctions can only be tangential to eigenvectors of $Df(x^*)^T$. This is why the following lemma is of importance.

\begin{lemma}\label{lem:DerivativeConjugated}
	Let $I$ induce a subsystem for $f$. Then for any $x \in X$
	\begin{equation}\label{eq:IntertwineDf}
	   Df(x)^T \cdot D\Pi_I^T = D\Pi_I^T \cdot Df_I(\Pi_I(x))^T \text{ and } D\Pi_I \cdot Df(x) = Df_I(\Pi_I(x)) D\Pi_I(x)
	\end{equation}
	where the derivate $D\Pi_I$ of $\Pi_I$ is the (constant) matrix with rows consisting of the standard basis vectors $(e_i)_{i \in I}$ and $\cdot$ denotes the matrix product. 
	In particular, if $w = (w_i)_{i \in I}$ is an eigenvector of $Df_I(\Pi_I(x))^T$ with eigenvalue $\lambda$ then so is $\overline{w}$ for $Df(x)^T$, for
	\begin{equation}\label{equ:ExtendEigenvector}
			\overline{w} := D\Pi_I^Tw = \begin{cases} w_k, & k \in I\\
											0, & \text{ else}
							  \end{cases}
	\end{equation}
	and if $v$ is an eigenvector of $Df(x)$ with eigenvalue $\lambda$ and $\Pi_I(v) \neq 0$ then $\Pi_I(v)$ is an eigenvector of $Df_I(\Pi_I(x))$ with eigenvalue $\lambda$.
\end{lemma}

\begin{proof} From the subsystem equation (\ref{equ:def:subsystem}) it follows from linearity of $\Pi_I$
\begin{equation}\label{eq:DifferentialConjugated}
	D\Pi_I \cdot Df = D(\Pi_I \circ f) = D(f_I \circ \Pi_I) = (Df_I \circ \Pi_I) \cdot D\Pi_I.
\end{equation}
We obtain (\ref{eq:IntertwineDf}) by taking the transpose in (\ref{eq:DifferentialConjugated}). If $w$ is an eigenvector of $Df_I(x)^T$ with eigenvalue $\lambda$ then we get
\begin{equation*}
    Df(x)^T \overline{w} = Df(x)^T  \cdot D\Pi_I^Tw = \Pi_I^T \cdot Df_I(x)^T w = \lambda \Pi_I^T w = \lambda \overline{w}.
\end{equation*}
If $v$ is an eigenvector of $Df(x)$ with eigenvalue $\lambda$ we have by (\ref{eq:IntertwineDf})
\begin{equation*}
    Df_I(\Pi_I(x)) \Pi_I(v) = Df_I(\Pi_I(x)) \cdot D\Pi_I v = D\Pi_I \cdot Df(x) v = \lambda D\Pi_I(v) = \lambda \Pi_I(v).
\end{equation*}
The condition $\Pi_I(v) \neq 0$ guarantees that $\Pi_I(v)$ is an eigenvector of $Df_I(\Pi_I(x))$.
\end{proof}

For uniqueness of eigenfunctions the concept of resonance is crucial. It is motivated by the property that for eigenfunctions $g_1,\ldots, g_k$ with eigenvalues $\lambda_1,\ldots,\lambda_k$ and $r_1,\ldots,r_k\in \N$ we have $\prod\limits_{i = 1}^k g_i^{r_i}$ is an eigenfunction with eigenvalue $\sum\limits_{i = 1}^k r_i\lambda_i$.

\begin{definition}[Resonance condition; \cite{MezicKoopmanOperatorSpectrum} p. 12]
	We say a matrix $A \in \mathbb{C}^{n\times n}$ with eigenvalues $\lambda_1,\ldots,\lambda_n$ (with algebraic multiplicity) is resonant if there exists some $i \in \{1,\ldots,n\}$ and $m_1,\ldots,m_n \in \N$ with $\sum\limits_{j = 1}^n m_j \geq 2$ such that
	\begin{equation*}
		\lambda_i = \sum\limits_{j \neq i}^n m_j \lambda_j.
	\end{equation*}
	We say $A$ is resonant of order $k$ if $(m_1,\ldots,m_n)$ can be chosen such that $\sum\limits_{j = 1}^n m_j \leq k$.
\end{definition}

Non-resonance (of order $k$), i.e. not being resonant (of order $k$), and regularity is what is needed to guarantee existence and uniqueness of principal eigenfunctions (\cite{Regularity}, \cite{Matthew}). When we apply this to subsystems we see that a principal eigenfunction uniquely corresponds to a principal eigenfunction for a subsystem and vice versa. In particular that $X$ factors implies that all principal eigenfunctions can be recovered from the subsystems. This is the statement of the following Theorem and Corollary \ref{cor:EigenfunctionSubsystemEigenfunction}.

\begin{theorem}\label{thm:PrincipalEigenvalues}
	Assume there exists a globally exponentially stable fixed point $x^*$ and let $I$ induce a subsystem. Assume $Df(x^*)^T$ is diagonalizable and, for simplicity, that all eigenvalues $\lambda_1,\ldots,\lambda_n$ have algebraic multiplicity one with corresponding eigenvectors $v_1,\ldots,v_n$. Assume that the eigenvalues are non-resonant of order $k$ with $k > \max\limits_{i,j} \frac{\mathrm{Re}(\lambda_i)}{\mathrm{Re}(\lambda_j)}$. Assume $f$ is $k$ times continuously differentiable. Then there exist $n$ uniquely determined principal eigenfunction of the whole system tangential to $v_1,\ldots,v_n$ and exactly $|I|$ pairwise distinct principal eigenfunctions for the subsystem, each tangential to one of the $\Pi_I(v_j)$ for some $j$, and they induce principal eigenfunctions for the whole system.
\end{theorem}

\begin{proof}
    By \cite{Matthew} Proposition 6, the assumptions guarantee the existence and uniqueness of $n$ principal eigenfunctions $g_1,\ldots,g_n$ tangential to $v_1,\ldots,v_n$. Next, we show that the assumption on non-resonance of $Df(x^*)$ carries over to $Df_I(\Pi_I(x^*))$ and we can use \cite{Matthew} Proposition 6 for the subsystem as well. By Lemma \ref{lem:DerivativeConjugated} we get that the spectrum of $Df_I(\Pi_I(x^*))$ is contained in the spectrum of $Df(x^*)$. Further, it also follows from Lemma \ref{lem:DerivativeConjugated} that the geometric multiplicity of each eigenvalue $\lambda$ for $Df_I(\Pi_I(x^*))$ is at most the geometric multiplicity of $\lambda$ for $Df(x^*)$, i.e. at most 1. Non-resonance (of order $k$) of $Df_I(\Pi_I(x^*))$ follows now from non-resonance (of order $k$) of $Df(x^*)$. For each basis of eigenvectors $w_1,\ldots,w_{|I|}$ of $Df_I(\Pi_I(x^*))^T$, we use \cite{Matthew} Proposition 6, to guarantee the existence and uniqueness $|I|$ many principal eigenfunctions $h_1,\ldots,h_{|I|}$ tangential $w_1,\ldots,w_{|I|}$. It remains to show that each of the $w_1,\ldots,w_{|I|}$ can be chosen to be of the form $w_j = \Pi_I(v_{i(j)})$ for $1\leq j \leq |I|$ and some $1\leq i(j) \leq n$. Lemma \ref{lem:DerivativeConjugated} states that for $j = 1,\ldots,|I|$, the vectors $D\Pi_I^T w_j$ are eigenvectors of $Df(x^*)^T$. From the assumption that all eigenvalues $\lambda$ of $Df(x^*)$ have algebraic (hence also geometric) multiplicity one, it follows that there exist unique $i(j) \in \{1,\ldots,n\}$ and $r_j \in \R$ such that $r_jD\Pi_I^T w_j = v_{i(j)}$. That means $r_jh_j$ is a principal eigenfunction tangential to $r_j w_j = \Pi_I(D\Pi_I^T (r_j w_j)) = \Pi_I(v_{i(j)})$. And Proposition \ref{propEigenvector} implies that $\tilde{g}_j := r_jh_j \circ \Pi_I$ is a principal eigenfunction ($\nabla \tilde{g}_j(x^*) = v_{i(j)} \neq 0$) of the whole system.
\end{proof}

The following corollary addresses the question whether we can find all the principle eigenfunctions by searching for them in subsystems. The answer is positive.

\begin{corollary}\label{cor:EigenfunctionSubsystemEigenfunction}
	Additionally to the assumptions of Theorem \ref{thm:PrincipalEigenvalues} assume that we have subsystems $(I_1,f_{I_1}),\ldots,(I_N,f_{I_N})$ such that $\bigcup\limits_{k = 1}^N I_k = \{1,\ldots,n\}$ then any principal eigenfunction for the whole system is already a principal eigenfunction for one of the subsystems.
\end{corollary}

\begin{proof}
Let $g$ be a principal eigenfunction of the whole system with $w := \nabla g(x^*)$. Then $w$ is an eigenvector of $Df(x^*)^T$ with eigenvalue $\lambda$. Hence $\lambda$ is also an eigenvalue of $Df(x^*)$. Let $v$ be its corresponding eigenvector. From $\bigcup\limits_{k = 1}^N I_k = \{1,\ldots,n\}$ it follows that for at least one $k \in \{1,\ldots,N\}$ we have $\Pi_{I_k}(v) \neq 0$. From Lemma \ref{lem:DerivativeConjugated} we get that $\lambda$ is an eigenvalue of $Df_{I_k}(\Pi_{I_k}(x^*))$ (with eigenvector $\Pi_{I_k}(v)$). Hence $\lambda$ is also an eigenvalue of $Df_{I_k}(\Pi_{I_k}(x^*))^T$. As in the proof of Theorem \ref{thm:PrincipalEigenvalues}, we see that there exists an eigenfunction $h$ with eigenvalue $\lambda$ for the subsystem induced by $I_k$ such that $\tilde{g} := h \circ \Pi_{I_k}$ is an eigenfunction (with eigenvalue $\lambda$) of the whole system. Because we assumed that the eigenvalues are simple, by scaling, we get $\nabla \tilde{g} (x^*) = \nabla g(x^*)$. Uniqueness of principal eigenfunctions implies $\tilde{g} = g$. That shows that $g$ is induced by a principal eigenfunction from a subsystem, namely $h$.

\end{proof}


Finding eigenfunctions for the subsystems is not answered by Theorem \ref{thm:PrincipalEigenvalues} and remains a general task (as for finding invariant measures). A partial answer to that question is given in \cite{MilanEigenfunction} and also Laplace averages can be used, where Proposition 6 and Remark 14 from \cite{Matthew} provides a condition under which the Laplace averages exist.

For applications a good choice of subsystems (including a factorization of $X$ in the case for invariant measures, see Theorem \ref{thm:GluedInvariantMeasure}) is essential. In case of finding (principal) eigenfunctions it is beneficial to find subsystems as small as possible because we have seen in Proposition \ref{propEigenvector} and Theorem \ref{thm:PrincipalEigenvalues} that (principal) eigenfunctions are inherited to larger systems. For such families of (nested) subsystems we refer to~\cite{SparsePaper}.

\subsection{Coordinate free formulations}\label{subsec:CoordinateFreeFormulation}

For the definition of subsystems we used the concept of coordinate indices $I$, which is coordinate-dependent. Here we describe a coordinate-free formulation using the concept of factor systems.

\begin{definition}\label{def:SubsystemManifold}
We call $(Y,\Pi)$ a smooth factor system of the dynamical system on a smooth (compact) manifold (with boundary) $X$ with a smooth semiflow $\varphi_t$, if $Y$ is a smooth manifold (with boundary) and $\Pi:X \rightarrow Y$ a surjective submersion, i.e. $\Pi$ is surjective and $D\Pi(x)$ has full rank for all $x \in X$, such that there exists a smooth semiflow $\varphi^{Y}_t$ on $Y$ with
\begin{equation}\label{equ:def:SubsystemManifold}
	\varphi^{Y}_t \circ \Pi = \Pi \circ \varphi_t \; \; \text{ for all } t\in \R_+.
\end{equation}
\end{definition}

The idea of a factor system is to transport a dynamical system into another (less complex or lower dimensional) dynamical system. By (\ref{equ:FlowCommuteSubsystem}) we see that subsystems are special cases of factor systems by setting $X = \R^n$, $Y = \R^{|I|}$, $\Pi = \Pi_I$ and $\varphi^Y = \varphi^{I}$.

\begin{remark}
We consider smooth factor systems instead of factor systems -- where no smoothness is assumed -- in order to rule out pathologies for $\Pi$ as for instance space filling curves. Since the dimension of the image of a smooth map can not exceed the dimension of the domain (by Sard's theorem, for instance) we see that for smooth factor systems the dimension of $Y$ is necessarily at most the dimension of $X$. Further smoothness is needed in order to formulate an analogue version of Theorem \ref{thm:PrincipalEigenvalues} were regularity played an essential role.
\end{remark}

Many -- but not all -- of our arguments can be carried out in a coordinate free formulation for factor systems. Even though Proposition \ref{propEigenvector} holds in the same way for factor systems, Theorem \ref{thm:GluedInvariantMeasure} is not generalized to factor systems in a straight forward way because in the proof of Theorem \ref{thm:GluedInvariantMeasure} we used Lemma 2.1 from \cite{AmbrosioOptimalTransport} which requires that the maps $\Pi$ are projections. On the other hand Theorem \ref{thm:PrincipalEigenvalues} is true for smooth systems $(X,(\varphi_t)_{t \in \R_+})$ with smooth factor systems, that means $X$ and $Y$ are smooth manifolds, the semiflows $\varphi_t$ and $\varphi_t^Y$ are smooth and $\Pi$ is a smooth submersion. 
In order to see that we have to generalize Lemma \ref{lem:DerivativeConjugated} to this case. Equation (\ref{equ:def:SubsystemManifold}) implies by taking the time derivative in $t = 0$ that
\begin{equation}\label{eq:IntertwinefManifoldSubsys}
	D\Pi(x) f(x) = f^Y(\Pi(x))
\end{equation}
where $f^Y$ denotes the vector field that induces the flow $\varphi_t^Y$. Taking the derivative in (\ref{eq:IntertwinefManifoldSubsys}) with respect to $x$, using $f(x^*) = 0$ and evaluating in $x^*$ respectively $\Pi(x^*)$ gives
\begin{equation*}
	D\Pi(x^*) Df(x^*) = D\Pi(x^*) Df(x^*) + D^2\Pi(x) f(x) = Df^Y(\Pi(x^*)) D\Pi(x^*)
\end{equation*}
which exactly gives the first statement in Lemma \ref{lem:DerivativeConjugated}. To conclude that $Df^Y$ is diagonalizable it is needed to assume the condition that $\Pi$ is a submersion around $x^*$, i.e. $D\Pi(x^*)$ has full rank. For the projections $\Pi_I$ this is automatically the case, hence we did not need this assumptions. If $D\Pi(x^*)$ has full rank then at least $\mathrm{dim}(Y)$ many eigenvectors of $Df(x^*)$ are mapped by $D\Pi(x^*)$ to non-zero (distinct) eigenvectors of $Df^Y(x^*)$ and we proceed by the same arguments as in the proof of Theorem \ref{thm:PrincipalEigenvalues}.

\section{Computational applications to dynamic mode decomposition and invariant measures}

In this section we show that a-priori knowledge of subsystems can be used to reduce computational complexity. We demonstrate this at the examples of dynamic mode decomposition and of computation of invariant measures. Compared to sparsity approaches from \cite[Chapter 9]{DMDBook}, we work explicitly with known sparsity patterns in the dynamics, i.e. the structure of the Koopman operator, instead of observed sparsity in the data, which assumes sparsity implicitly. Sparsity in the dynamics naturally leads to sparsity in the data. But there are two reasons why the a-priori knowledge about sparse patterns is useful: First, it allows to enforce specific structures on the approximating object, here the approximation of the Koopman operator, based on the sparsity structures of the dynamics, without cost of accuracy. Second, the correct sparse structures might not be so easily detected just from data.

From Proposition \ref{propEigenvector} we obtain that the composition operator $V_I$ intertwines the systems. This structure can be enforced on the finite dimensional approximation of the Koopman operator a-priori and leads to lower dimensional systems and hence reduces computational complexity. We demonstrate that for dynamic mode decomposition at the example of coupled Duffing equations of the form (\ref{equ:Sparsex0x1x2}) with
\begin{align}\label{equ:Duffing12}
	\dot{x}^1_1 & = x^1_2\notag\\
	\dot{x}^1_2 & = -\delta x^1_2 - x^1_1 (\beta + \alpha)(2x^1_1)^2
\end{align}
for $\delta = 0.5$, $\beta = -1$ and $\alpha = 1$ and for $x^2,x^3$ the dynamics is for $i = 2,3$ for $\gamma_1 = 1$ and $\gamma_2 = 2$
\begin{align}\label{equ:Duffing3456}
	\dot{x}^i_1 & = x^i_2 \notag\\
	\dot{x}^i_2 & = -\delta x^i_2 - x^i_1 (\beta + \alpha)(2x^i_1)^2 + \gamma_i x_1^1
\end{align}
Other examples that have non-trivial subsystems include Dubins Car, 6D Acrobatic Quadrotor (see \cite{Chen}) or radial distribution power networks \cite{RadialDistributionNetworks}.

\subsection{Sparse dynamic mode decomposition}

In this section we describe how extended dynamic mode decomposition (EDMD) (\cite{DMDSchmid}) can benefit from exploiting subsystems, i.e. in this context utilising knowledge of subsystems. Theorems \ref{thm:GluedInvariantMeasure} and \ref{thm:PrincipalEigenvalues} indicate that this allows to still capture (some) important spectral properties of the dynamical system.

We propose a sparse EDMD, illustrated at the example of the coupled Duffing equations (\ref{equ:Duffing12}) and (\ref{equ:Duffing3456}) but it can be extended in the same manner to general systems with more complicated subsystem.

The idea is simple: Instead of applying the EDMD to the whole system we use EDMD separately for each of the separate subsystems.

\begin{figure}[!h]
\begin{picture}(250,220)
\put(10,60){\includegraphics[width=78mm]{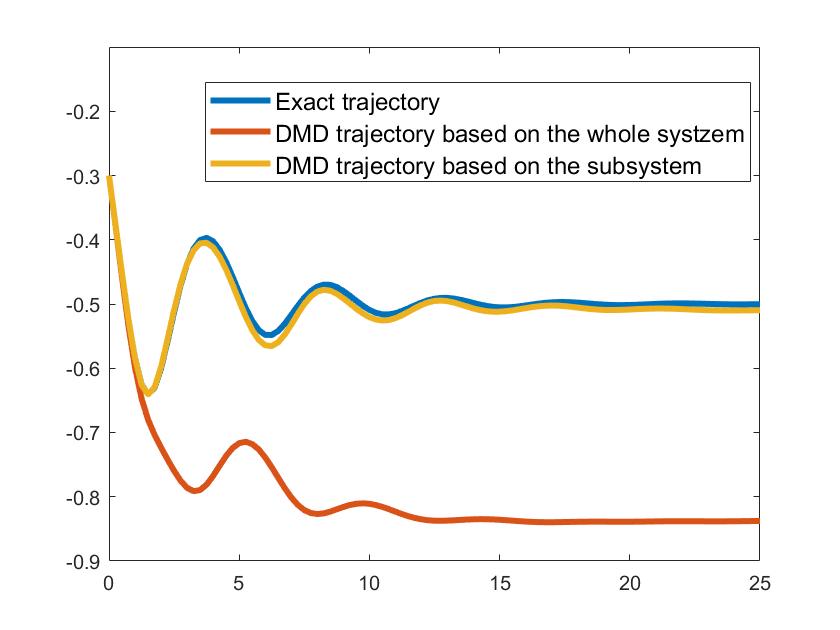}}
\put(220,60){\includegraphics[width=78mm]{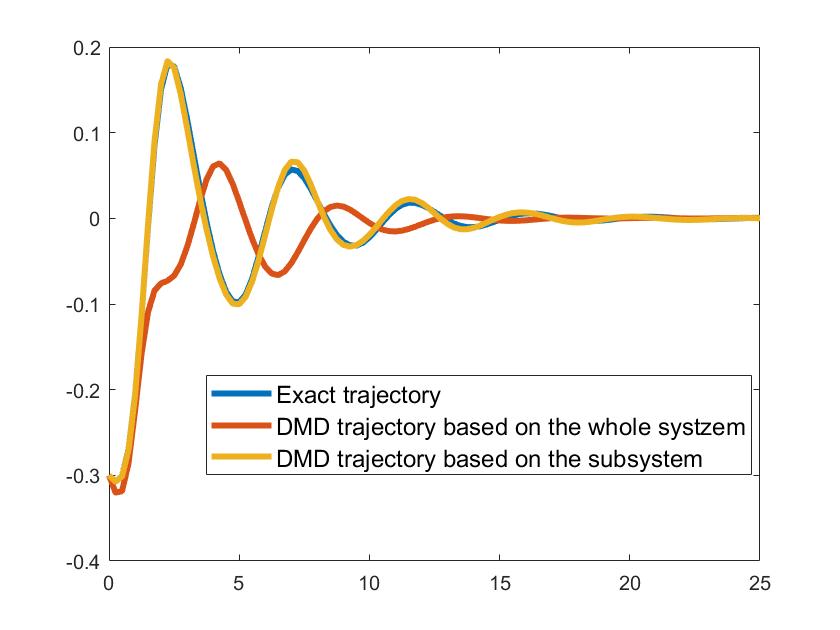}}
\end{picture}
\vspace{-20mm}
\caption{\footnotesize{DMD approximations based on the whole system and on the subsystem induced by $I = \{1,2\}$ for (\ref{equ:Duffing12}) with initial value $x_0 = (-0.3,-0.3,0.7,0.5,0.3,0.2)$, trained on the same data. Left: DMD approximation for $x_1^1$, right: DMD approximation for $x^1_2$. For the DMD for the whole system 1000 randomly generated radial basis functions were used while 350 radial basis functions were sufficient for the subsystem.}}
\label{figDufffing12}
\end{figure}

\begin{figure}[!h]
\begin{picture}(250,220)
\put(10,60){\includegraphics[width=78mm]{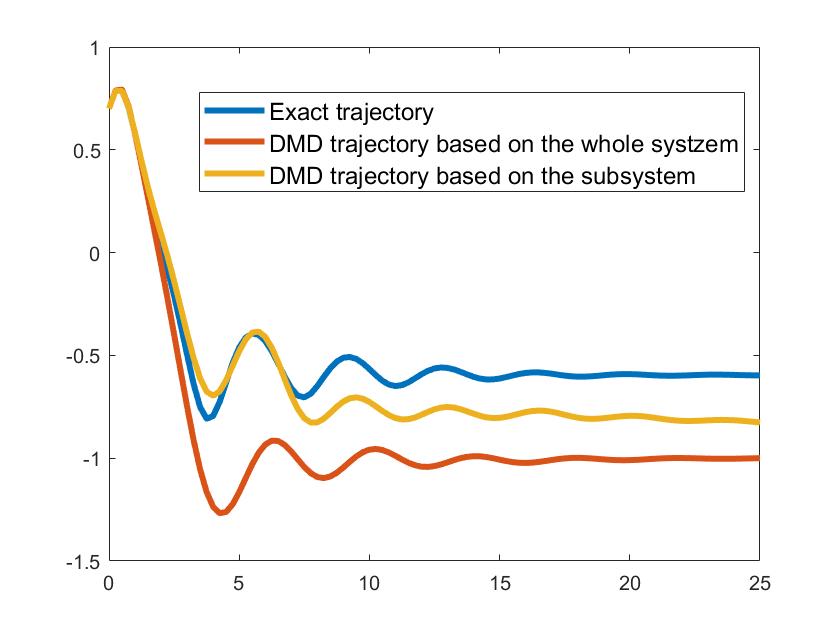}}
\put(220,60){\includegraphics[width=78mm]{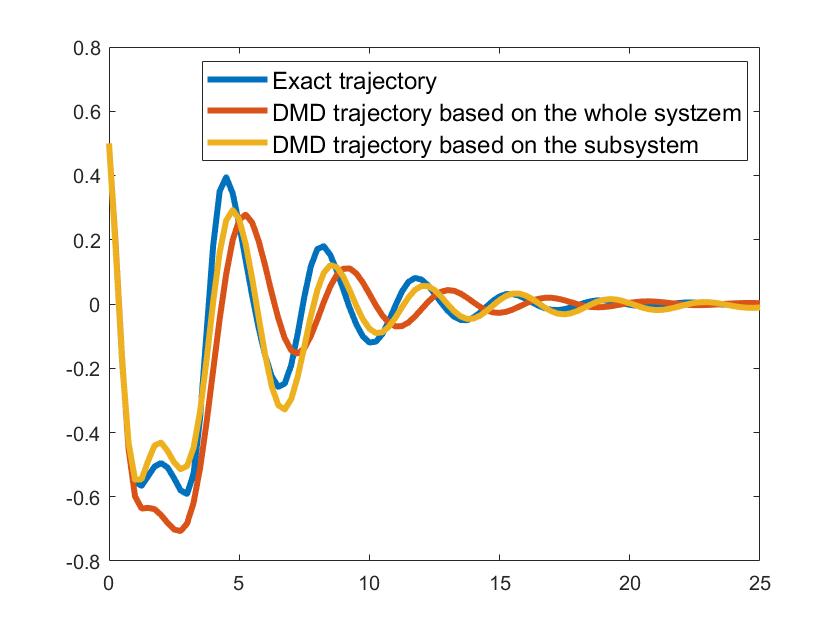}}
\end{picture}
\vspace{-20mm}
\caption{\footnotesize{DMD approximations based on the whole system and on the subsystem induced by $I = \{1,2,3,4\}$ for (\ref{equ:Duffing12}) with initial value $x_0 = (-0.3,-0.3,0.7,0.5,0.3,0.2)$, trained on the same data with randomly generated 1000 radial basis functions. Left: DMD approximation for $x_2^1$, right: DMD approximation for $x^2_2$.}}
\label{figDufffing34}
\end{figure}

\begin{figure}[h]
\begin{picture}(250,220)
\put(10,60){\includegraphics[width=78mm]{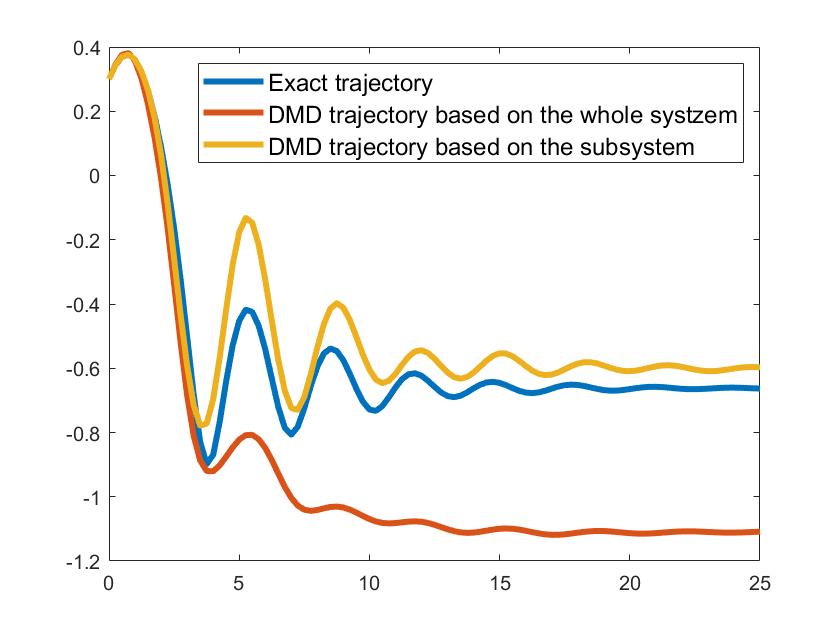}}
\put(220,60){\includegraphics[width=78mm]{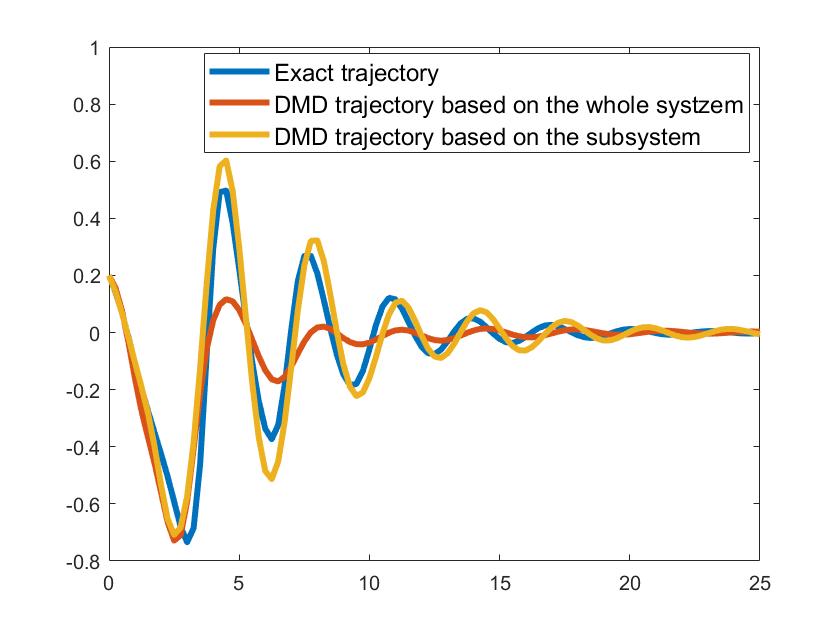}}
\end{picture}
\vspace{-20mm}
\caption{\footnotesize{DMD approximations based on the whole system and on the subsystem induced by $I = \{1,2,5,6\}$ for (\ref{equ:Duffing12}) with initial value $x_0 = (-0.3,-0.3,0.7,0.5,0.3,0.2)$, trained on the same data with randomly generated 1000 radial basis functions. Left: DMD approximation for $x_3^1$, right: DMD approximation for $x^3_2$.}}
\label{figDufffing56}
\end{figure}

Depending on the choice of dictionary, i.e. functions $\mathbf{\Psi}= (\psi_1,\ldots,\psi_l)$ on which we apply EDMD, this comes with the following advantage
\begin{enumerate}
    \item In case the number $l$ of dictionary functions is fixed: On a lower dimensional space the same number of dictionary functions allows a better resolution of the geometry the space. For example if we choose a dictionary with 1000 radial basis functions (as used to obtain Figures \ref{figDufffing12}, \ref{figDufffing34} and \ref{figDufffing56}) for the whole state space $\R^6$ as well as for $\R^2$ and respectively $\R^4$ we get a better geometric description of the spaces $\R^2$ and $\R^4$ compared to $\R^6$, allowing a finer EDMD approximation; see Figures \ref{figDufffing12}, \ref{figDufffing34} and \ref{figDufffing56}.
    \item In case the number $l$ of dictionary functions depends on the dimension of the space: Typically $l$ is larger for larger dimensions. This, for instance, is the case when the dictionary consists of (trigonometric) polynomials up to a certain degree. The dimension of the space of polynomials of degree up to $d$ in dimension $n$ is given by $\binom{n+d}{n}$ and hence grows combinatorial in the space dimension $n$. This relates to the curse of dimensionality and underlines the benefitial impact of lowering the dimension of the space.
\end{enumerate}
In the following we will describe in more detail how we propose to exploit sparsity for EDMD at the example of the coupled Duffing equations (\ref{equ:Duffing12}) and (\ref{equ:Duffing3456}). For a simpler notation we denote the pairs $(x^i_1,x^i_2) \in \R^2$ by $x^i \in \R^2$. Then we use the subsystems induced by $I_1 = \{1\}$, $I_2 = \{1,2\}$ and $I_3 = \{1,3\}$ representing the subsystems on the sates $x^1 = (x^1_1,x^1_2)$ and $(x^1,x^i) = (x^1_1,x^1_2,x^i_1,x^i_2)$ for $i = 2,3$.

For EDMD we are given a data set of snapshots $x(k) =(x^1(k),x^2(k),x^3(k))^T \in \R^6$ for $k = 0,\ldots,m$ and corresponding evolutions $y(k) = f(x(k))$ for $k = 0,\ldots,m$. Further, we have a dictionary $\mathbf{\Psi} = (\psi_1,\ldots,\psi_l)$ of functions on $\R^6$. As mentioned before in our example for the coupled Duffing equation we choose $\mathbf{\Psi}$ to be a dictionary of $l = 1000$ Gaussian radial basis functions with randomly generated centers in $\R^6$. The EDMD for the whole system has the following form
\begin{equation}\label{DMDWhole}
	\mathbf{K} \in \mathrm{arg}\min\limits_{A \in \R^{l \times l}} \sum\limits_{k = 0}^m \left\|\mathbf{\Psi}(y(k)) - A\mathbf{\Psi}(x(k))\right\|^2
\end{equation}
and generates a finite dimensional approximation $\mathbf{K}$ of the Koopman operator.


Our sparse EDMD for the coupled Duffing equations (\ref{equ:Duffing12}), (\ref{equ:Duffing3456}) is stated in Algorithm \ref{alg:SparseDMD}.
\begin{algorithm}\label{Alg:SparseDMD}
\caption{Sparse dynamic mode decomposition}
\label{alg:SparseDMD}
\begin{algorithmic}[1]
\BState \emph{Input:} Snapshots $x(k)=(x^1(k),x^2(k),x^3(k))^T$ for $k = 0,\ldots, m$
\BState \emph{Split the snapshots into snapshots of the subsystems:}
	Split $x(k)$ into $x^1(k)$ and $(x^1(k),x^2(k))$ respectively $(x^1(k),x^3(k))$ for $k = 0,\ldots,m$.
\BState \emph{Choose dictionaries:}
	Let $\mathbf{\Psi}^1 = (\psi^1_1,\ldots,\psi^1_l)$, $\mathbf{\Psi}^{1,i} = (\psi^{1,i}_1,\ldots,\psi^{1,i}_l)$, for $i = 2,3$ for $l = 1000$, be Gaussian radial basis functions with randomly generated centers in $\R^2$ respectively $\R^4$.
\BState \emph{Compute approximation matrices for the Koopman operators for the subsystems:} Compute the approximating matrices $\mathbf{K}^1,\mathbf{K}^{1,2},\mathbf{K}^{1,3}$ for the Koopman operator based on the corresponding snapshots and dictionaries, i.e.
    \begin{equation}\label{DMD1}
	    \mathbf{K}^1 \in \mathrm{arg}\min\limits_{A\in \R^{l \times l}} \sum\limits_{k = 0}^m \left\|\mathbf{\Psi}^1(y^1(k)) - A\mathbf{\Psi}^1(x^1(k))\right\|^2
    \end{equation}
    for the first subsystem, and
    \begin{equation}\label{DMD2}
	    \mathbf{K}^{1,i} \in \mathrm{arg}\min\limits_{A\in \R^{l \times l}} \sum\limits_{k = 0}^m \left\|\mathbf{\Psi}^{1,i}(y^1(k),y^i(k)) - A\mathbf{\Psi}^{1,i}(x^1(k),x^i(k)) \right\|^2
    \end{equation}
    for the second, respectively third subsystem.
\BState \emph{Return} $\mathbf{K}^1, \mathbf{K}^{1,2}, \mathbf{K}^{1,3}$
\end{algorithmic}
\end{algorithm}

As mentioned, the advantage of solving the separate optimization problems lies in either having a lower dimensional optimization problem to solve or a better resolution of the lower dimensional state space. The impact of these advantages gets more and more important when the state space gets higher dimensional, as for example in fluid dynamics (\cite{IgorFluid}), or when the size of the dictionary in the EDMD depends on the state space dimension.

The approximations of the Koopman operator for the subsystems can be used for state estimation or computation of (principal) eigenfunctions. This is based on (\ref{equ:FlowCommuteSubsystem}), Proposition \ref{propEigenvector} and Theorem \ref{thm:PrincipalEigenvalues}. We illustrate state estimation in Figures \ref{figDufffing12}, \ref{figDufffing34} and \ref{figDufffing56}, which display numerical examples of our approach for the coupled Duffing equations (\ref{equ:Duffing12}) and (\ref{equ:Duffing3456}). We use 500 sample trajectories sampled by step size $0.25$ for $25$ time steps and a vector $g$ of $1000$ Gaussian radial basis functions with randomly generated centers. Figures \ref{figDufffing12}, \ref{figDufffing34} and \ref{figDufffing56} display a comparison of the state estimation via EDMD and sparse EDMD. We took the initial value $x_0  = (-0.3,-0.3,0.7,0.5,0.3,0.2)$ that has not occurred in the snapshots and use $\mathbf{K}$ (from the non-sparse EDMD (\ref{DMDWhole})) to estimate all the states $x_0(t) = (x^1_0(t),x_0^2(t),x_0^3(t))$ together for $25$ time units. In the sparse EDMD we estimate the sate $x_0^1(t)$, using only $\mathbf{K}^1$ and the dictionary $\mathbf{\Psi}^1$, the state $x_0^2(t)$, using $\mathbf{K}^{1,2}$ obtained from the dictionary $\mathbf{\Psi}^{1,2}$ and $x_0^3(t)$, using $\mathbf{K}^{1,3}$ obtained from the dictionary $\mathbf{\Psi}^{1,3}$, for $25$ time steps.


\subsection{Sparse computation of invariant measures}

In this section we propose a sparse computation of invariant measures for the approach from \cite{VictorInvariantMeasure} and \cite{MilanInvariantMeasure} for polynomial dynamical systems based on convex optimization. Before stating the approach we want to shortly present the underlying idea. For this purpose it is easier to consider discrete dynamical systems, i.e. $x_{k+1} = f(x_k)$ for some continuous function $f:\R^n \rightarrow \R^n$. As mentioned in Remark \ref{RemarkDiscrete} sparse properties of the Koopman operator hold in an analogous way to continuous time systems for discrete time systems as well. A measure $\mu$ on $X$ is invariant if and only if we have
\begin{equation}\label{equ:InvMeasureDense}
	\langle g, \mu \rangle = \int\limits_X g \; d\mu = \int\limits_X g \circ f\; d\mu = \langle g \circ f, \mu \rangle
\end{equation}
for all $g$ in (a dense subset of) $\C(X)$. To follow \cite{MilanInvariantMeasure} we formulate the following linear optimization problem for extremal invariant probability measures
\begin{equation}\label{equ:OptProbInvMeasure}
\begin{tabular}{lllll}
		$p^*$ & $:=$ & $\min $ & $\langle G, \mu \rangle = \int\limits_X G \; d\mu$&\\
		& & s.t. & $ \mu \in \mathcal{M}(X)_+$ &\\
			& & & $ \mu(X) = 1$\\
			& & & $\langle g,\mu\rangle = \langle g \circ f,\mu\rangle$ for all $g \in \C(X)$ &
	\end{tabular}
\end{equation}
Where $G$ represents a cost and can be used to identify specific invariant measures that we are interested in since invariant measures are not unique in general. Now assume $I_1,\ldots,I_N$ induce subsystems and the cost $G$ is adapted to the sparse structure -- that means $G$ can be written as $G = \sum\limits_{k = 1}^N G_k$ for $G_k \in \C(\Pi_{I_k}(X))$ for $k = 1,\ldots,N$. We formulate corresponding sparse linear programming problems (\ref{equ:OptProbInvMeasure1}) and (\ref{equ:OptProbInvMeasure2}) given by
	\begin{equation}\label{equ:OptProbInvMeasure1}
\begin{tabular}{lllll}
		$p_1^*$ & $:=$ & $\min$ & $\int\limits_X G \; d\mu$&\\
		& & s.t. & $ \mu \in \mathcal{M}(X)_+$ &\\
			& & & $ \mu(X) = 1$\\
			& & & $ \int\limits_{X} h \circ f_{I_k} \; d\mu = \int\limits_{X} h \; d\mu$ & for all $h \in \C(\Pi_{I_k}(X))$ for $k = 1,\ldots,N$
	\end{tabular}
\end{equation}
and
\begin{equation}\label{equ:OptProbInvMeasure2}
\begin{tabular}{lllll}
		$p_2^*$ & $:= $ & $\min$ & $\sum\limits_{k = 1}^N  \langle G_k, \mu_k\rangle $&\\
		& & s.t. & $ \mu_k \in \mathcal{M}(\Pi_{I_k}(X))_+$ & $k = 1,\ldots,N$\\
			& & & $ \langle \mathbf{1},\mu_k\rangle = 1$ & $k = 1,\ldots,N$\\
			& & & $ \langle h \circ f_{I_k},\mu_k\rangle = \langle h, \mu_k\rangle$ & $\forall h \in \C(\Pi_{I_k}(X))$ for $k = 1,\ldots,N$\\
			& & & $ \langle h, \mu_k\rangle = \langle h, \mu_l\rangle$ & $\forall h \in \C(\Pi_{I_k \cap I_l}(X))$ for $k,l = 1,\ldots,N$.
	\end{tabular}
\end{equation}

The linear programming problem (\ref{equ:OptProbInvMeasure1}) is clearly a relaxation of (\ref{equ:OptProbInvMeasure}) since invariance is only required for the marginals of corresponding to the subsystems. The linear programming problem (\ref{equ:OptProbInvMeasure2}) is a reduction to a vector of measures on lower dimensional spaces, and is what we recommend for practical computations, while (\ref{equ:OptProbInvMeasure1}) allows more degrees of freedom. 

\begin{proposition}\label{prop:EquivalentExtremalMeasures}
	Let $(I_1,f_{I_1}),\ldots, (I_N,f_{I_N})$ be subsystems such that $X$ factors with respect to them. Let $G \in \C(X)$ such that $G$ can be written as $G = \sum\limits_{k = 1}^N G_k \circ \Pi_{I_k}$ for $G_k \in \C(\Pi_{I_k}(X))$ for $k = 1,\ldots,N$. Then for (\ref{equ:OptProbInvMeasure}) we have $p^* = p_1^* = p_2^*$ from (\ref{equ:OptProbInvMeasure1}) and (\ref{equ:OptProbInvMeasure2}). Further, there exists an invariant measure $\mu$ for the whole system with $p^* = \langle G, \mu \rangle$ such that $\mu$ is optimal for (\ref{equ:OptProbInvMeasure1}) and $\mu_k := (\Pi_{I_k})_\#\mu$ for $k = 1,\ldots,N$ form an optimal (feasible) point for (\ref{equ:OptProbInvMeasure2}).
\end{proposition}

\begin{proof} We have already noted that (\ref{equ:OptProbInvMeasure1}) is a relaxation of (\ref{equ:OptProbInvMeasure}), i.e. $p^*_1 \leq p^*$. We also see that each feasible point $\mu$ for (\ref{equ:OptProbInvMeasure1}) induces a feasible point $(\mu_1,\ldots,\mu_N)$ for (\ref{equ:OptProbInvMeasure2}) by $\mu_k = (\Pi_{I_k})_\# \mu$. It follows
\begin{equation}\label{eq:p2p1p}
    p_2^* \leq p_1^* \leq p^*
\end{equation}
and it remains to show $p_2^* \geq p^*$. Therefore, let $(\mu_1,\ldots,\mu_N)$ be feasible for (\ref{equ:OptProbInvMeasure2}). Theorem \ref{thm:GluedInvariantMeasure} implies that we can find an invariant measure $\mu$ for the whole system with corresponding marginals $(\mu_1,\ldots,\mu_N)$. Since $G$ is sparse we get
\begin{eqnarray*}
	p^* & \leq & \langle G, \mu \rangle = \sum\limits_{k = 1}^N \langle G_k \circ \Pi_{I_k}, \mu\rangle = \sum\limits_{k = 1}^N \langle G_k, (\Pi_{I_k})_\# \mu\rangle = \sum\limits_{k = 1}^N \langle G_k, \mu_k\rangle.
\end{eqnarray*}
Since $(\mu_1,\ldots,\mu_N)$ was an arbitrary feasible point for (\ref{equ:OptProbInvMeasure2}) we get $p^* \leq p_2^*$. In fact, we have shown that an optimal point $(\mu_1,\ldots,\mu_N)$ for (\ref{equ:OptProbInvMeasure2}) is induced by an optimal point $\mu$ for (\ref{equ:OptProbInvMeasure}). The same measure is also optimal for (\ref{equ:OptProbInvMeasure}).
\end{proof}

\begin{remark}
	If any of the subsystems has a unique invariant measure this adds even more sparse structure to the optimization problem. 
\end{remark}

\subsubsection{Sparse moment-sum-of-squares formulation for extremal invariant measures}

In this section we describe how to use a sparse version of moment-sum-of-squares (moment-SOS) techniques from \cite{MilanInvariantMeasure}. We follow the notion from \cite{MilanInvariantMeasure} and introduce the necessary objects.

Before incorporating sparse structure we start by describing the moment-SOS hierarchy, that is, a semidefinite-programming representable approximation of $\mathcal{M}(Y)_+$, the space of non-negative measures, on a basic semialgebraic set $Y$.

Let $Y$ be a basic semialgebraic set, that is, $Y$ is given by
\begin{equation}\label{eq:YBasicSemialgebraic}
    Y = \{x \in \R^n: g_k(x) \geq 0, k = 1,\ldots,m\}
\end{equation}
for $m \in \N$ and polynomials $g_1,\ldots,g_m$.

Using so-called moment sequences, we define outer approximation $\mathcal{M}_d(Y)$ of $\mathcal{M}(Y)_+$ for even degree $d \in \N$ by
\begin{equation}\label{eq:M_d}
    \mathcal{M}_d(Y) := \{\mathbf{y} \in \R^{\N^n} : M_d(g_i\mathbf{y}) \succeq 0, i = 0,1,\ldots,m\}
\end{equation}
where $\R^{\N^n}$ denotes the space of sequences $\mathbf{y} = (y_\alpha)_{\alpha \in \N^n}$ and $y_\alpha \in \R$ for all $\alpha \in \N^n$, the function $g_0$ is the constant one polynomial with $g(x) = 1$ for all $x \in \R^n$, $\succeq$ denotes positive semidefiniteness of a matrix and $M_d(g_i \mathbf{y})$ the so-called moment-localizing matrices, which will be defined below.

To define the moment-localizing matrices $M_d(g_i\mathbf{y})$ we need the notion of Riesz functionals $\ell_\mathbf{y}^d:\R[x]_d \rightarrow \R$ for a given $\mathbf{y} \in \R^{\N^n}$ defined by
\begin{equation}\label{eq:RieszFunctional}
    \ell_\mathbf{y}^d\left(\sum\limits_{|\alpha| \leq d} p_\alpha x^\alpha \right) := \sum\limits_{|\alpha| \leq d} p_\alpha y_\alpha.
\end{equation}
Denoting $v_d(x):= (x^\alpha)_{|\alpha|\leq d}$ the sequence of monomials up to degree $d$ we define $M_d(g_i\mathbf{y})$ by
\begin{equation}\label{eq:DefLocalizingMatrix}
    M_d(g_i\mathbf{y}) := \ell_\mathbf{y}^d(g_iv_{d_i}v_{d_i}^T)
\end{equation}
where $d_i := \lfloor{\frac{d-\mathrm{deg}g_i}{2}}\rfloor$.

The relation between $\mathcal{M}(Y)_d$ and non-negative measures is given by integrating the measures against polynomials. That is, for a given measure $\mu \in \mathcal{M}(Y)_+$ (with $Y$ compact) we can define $\mathbf{y} = (y_\alpha)_{\alpha \in \N^n}$ by
\begin{equation}\label{eq:MomentY}
    y_\alpha:= \int\limits_Y x^\alpha \; d\mu.
\end{equation}
Then $\ell^d_\mathbf{y} (p) = \int\limits_Y p \; d\mu$ for all $p \in \R[x]_d$. From the condition $g_i(x) \geq 0$ on $Y$ it follows directly that $\ell^d_\mathbf{y}(g_ip^2) = \int\limits_Y g_i p^2 \; d\mu \geq 0$ for all $p \in \R[x]_{d_i}$ for $d_i := \lfloor{\frac{d-\mathrm{deg}g_i}{2}}\rfloor$. This shows that the positive semidefiniteness condition on $M_d(g_i\mathbf{y})$ is necessary for $\mathbf{y}$ to be a moment sequence of a non-negative measure $\mu \in \mathcal{M}(Y)_+$ as in (\ref{eq:MomentY}); conversely the moment version of Putinar's Positivstellensatz \cite{Putinar} states that if all the moment-localizing matrices $M_d(g_i\mathbf{y})$ are positive semidefinite (and one of the $g_i$ is of the form $g_i = R-\|x\|^2$ for some $R>0$) then $\mathbf{y} = (y_\alpha)_{\alpha \in \N^n}$ is given by (\ref{eq:MomentY}) for a unique measure $\mu \in \mathcal{M}(Y)_+$.

Let us turn back to dynamical systems. Let the constraint set $X$ be compact and semialgebraic (and without loss of generality we assume one of the $g_i$ in the semialgebraic representation of $X$ is of the form $g_i = R-\|x\|^2$). Further, we consider polynomial vector fields $f$. That allows us to reformulate the following optimization problem (from (\ref{equ:OptProbInvMeasure}))
\begin{equation*}
\begin{tabular}{lllll}
		$p^*$ & $:=$ & $\min $ & $\langle G, \mu \rangle = \int\limits_X G \; d\mu$&\\
		& & s.t. & $ \mu \in \mathcal{M}(X)_+$ &\\
			& & & $ \mu(X) = 1$\\
			& & & $\langle g,\mu\rangle = \langle g \circ f,\mu\rangle$ for all $g \in \C(X)$ &
	\end{tabular}
\end{equation*}
in terms of moment sequences by
\begin{equation}\label{equ:OptProbInvMeasureMoments}
	\begin{tabular}{llll}
	$p^*$ & $=$ & $\min$ & $\ell_\mathbf{y}^{\mathrm{deg}(G)}(G)$\\
		& s.t. & $ \mathbf{y} = (y_\alpha)_{\alpha \in \N^n}$ &\\
		   & & $(y_\alpha)_{|\alpha|\leq d} \in M_d(X)$ & for all $d \in \N$\\
		   & & $ y_0 = 1$&\\
	 	   & & $ \ell^d_\mathbf{y}(x^\alpha \circ f) = y_\alpha$ & for all $d \in \N$, for all $\alpha$ s.t. $\mathrm{deg}(x^\alpha \circ f) \leq d$.
	\end{tabular}
\end{equation}

By truncating the degree $d \in \N$ the following moment-SOS hierarchy of optimization problems is obtained (for $d\geq \mathrm{deg}(G)$)
\begin{equation}\label{equ:OptProbInvMeasureMomentsTruncated}
	\begin{tabular}{llll}
	$p_d^*$ & $:=$ & $\min$ & $\ell_\mathbf{y}^d(G)$\\
		& s.t. & $ \mathbf{y} = (y_\alpha)_{|\alpha|\leq d}$ &\\
		   & & $(y_\alpha)_{|\alpha|\leq d} \in M_d(X)$ & \\
		   & & $ y_0 = 1$&\\
	 	   & & $ \ell^d_\mathbf{y}(x^\alpha \circ f) = y_\alpha$ & for all $\alpha$ s.t. $\mathrm{deg}(x^\alpha \circ f) \leq d$.
	\end{tabular}
\end{equation}
Because any $\mathbf{y}$ feasible for (\ref{equ:OptProbInvMeasure}) induces a feasible point for (\ref{equ:OptProbInvMeasureMomentsTruncated}) by truncating we get $p_d^* \leq p^*$. By \cite{MilanInvariantMeasure} we also have $p_d^* \nearrow p^*$ as $d\rightarrow \infty$ by the aide of the moment version of Putinar's Positivstellensatz \cite{Putinar}.

The striking advantage of (\ref{equ:OptProbInvMeasureMomentsTruncated}) is that it is a finite dimensional semidefinite program, in particular convex, for which fast solvers exist and it provides a converging approximation of $p^*$ as $d\rightarrow \infty$. A drawback of this moment-SOS approach is that the number of variables $y_\alpha$ grows combinatorial in the dimension $n$, namely by $\binom{n+d}{n}$.

In order to reduce dimension we now exploit sparse structure in form of subsystems. Therefore, we combine the above moment-hierarchy with the result from Proposition \ref{prop:EquivalentExtremalMeasures} to consider moment-SOS hierarchies on several but lower dimensional spaces (namely on $\Pi_{I_k}(X)$ induced by subsystems $(I_1,f_{I_1}),\ldots,(I_N,f_{I_N})$).

Therefore, we assume a compatible structure of $X$, particularly a combination of algebraic and sparse structure.

\textbf{Assumption:} The set $X$ is compact and semialgebraic and factors with respect to the subsystems $(I_1,f_{I_1}),\ldots,(I_N,f_N)$.

Let $X$ be given by $X = \{x \in \R^n: g_1(x)\geq 0,\ldots,g_m(x)\geq 0\}$ for $m \in \N$ and polynomials $g_1,\ldots,g_m \in \R[x_1,\ldots,x_n]$. Sparse descriptions of $X$ are typically of the form
\begin{equation}\label{eq:XSparse}
    X = \{x \in \R^n: g_j(\Pi_{L_j}(x)) \geq 0 \text{ for } j = 1,\ldots,m\}
\end{equation}
where for $j = 1,\ldots,m$ we have $L_j \subset \{1,\ldots,n\}$ and $g_j$ are polynomials that depend only on the variables indexed by the sets $L_j$. Note that even if the sets $L_j$ induce subsystem it does not follow that $X$ factors with respect to the subsystems induced by the sets $L_j$ (see the short discussion after Definition \ref{def:Factorization}).

That $X$ factors with respect to $(I_1,f_{I_1}),\ldots,(I_N,f_N)$ means
\begin{equation*}
    X = \{x \in \R^n: \Pi_{I_k}(x) \in \Pi_{I_k}(X) \text{ for } k = 1,\ldots,N\}.
\end{equation*}
By the projection theorem from real algebraic geometry, the sets $\Pi_{I_k}(X)$ are semialgebraic for $k = 1,\ldots,N$. Therefore, we assume to have an explicit semialgebraic representation of $\Pi_{I_k}(X)$ for $k = 1,\ldots,N$ by
\begin{equation}\label{equ:PIkXSemialgebraic}
	\Pi_{I_k}(X) = \{x \in \R^{|I_k|} : g_0^{(k)}(x) \geq 0,\ldots,g_{l_k}^{(k)}(x) \geq 0\}
\end{equation}
for $k = 1,\ldots,N$ for some $l_k \in \N$ and polynomials $g_1^{(k)},\ldots,g_{l_k}^{(k)}$ on $\R^{|I_k|}$ where for each $k$ one of the $g^{(k)}_l$ is of the form $g^{(k)}(x) = R_k-\|x\|^2$ for some $R_k > 0$ and $g_0^{(k)} = 1$.

By virtue of Proposition \ref{prop:EquivalentExtremalMeasures} we can state the sparse optimization problem truncated at degree $d\in \N$. We can restate (\ref{equ:OptProbInvMeasureMomentsTruncated}), the moment optimization problem truncated at degree $d \in \N$, by
\begin{equation}\label{equ:OptLinearFunctionalTruncated}
\begin{tabular}{llllr}
		$s^*_d$ & $:= $ & $\min$ & $\sum\limits_{k = 1}^N  l_{\mathbf{y}^{(k)}}^d(G_k)$&\\
		& & s.t. & $ \mathbf{y}^{(k)} = (y^{(k)}_\alpha)_{\alpha}$ for $\alpha \in \N^{|I_k|}$, $|\alpha| \leq d$ & $k = 1,\ldots,N$\\
			& & & $  y^{(k)}_0 = 1$ & $k = 1,\ldots,N$\\
			& & & $ \ell^d_{\mathbf{y}^{(k)}}(x^\alpha \circ f) = y^{(k)}_\alpha $&  for $k = 1,\ldots,N$\\
			& & & $ \mathbf{y}^{(k)}_\alpha = \mathbf{y}^{(l)}_\alpha$ for $\alpha \in \N^{I_k \cap I_l}$, $|\alpha| \leq d$ & for $k,l = 1,\ldots,N$\\
			& & & $M_d(g_i^{(k)} \mathbf{y}^{(k)}) \succeq 0$ & $k = 1,\ldots,N, i = 0,\ldots,l_k$
	\end{tabular}
\end{equation}

Compared to the full semidefinite program (\ref{equ:OptProbInvMeasureMomentsTruncated}) with in total $\binom{n + \frac{d}{2}}{n}$ variables the sparse program has $\sum\limits_{k = 1}^N \binom{|I_k| + \frac{d}{2}}{|I_k|}$ many variables; which is significantly less than $\binom{n + \frac{d}{2}}{n}$ if $|I_k|$, the dimension of $\R^{|I_k|}$, is essentially smaller than $n$.

The following proposition states that as $d\rightarrow \infty$ the optimal values $s^*_d$ of the sparse problems converge to the optimal value $p^*$ of the infinite dimensional problem.

\begin{proposition}
	Let $X$ be positively invariant, compact semialgebraic and factor with respect to the subsystems $(I_1,f_{I_1}),\ldots, (I_N,f_{I_N})$. Let $f$ be polynomial. Without loss of generality we assume that for each $k = 1,\ldots,N$ one of the $g_i^{(k)}$ from (\ref{equ:PIkXSemialgebraic}) is of the form $g_i^{(k)} = r_k - \|\cdot\|_2^2$ for some $r_k >0$. Let $G = \sum\limits_{k = 1}^N G_k \circ \Pi_{I_k}$ for $G_k \in \R[(x_i)_{i \in I_k}]$ for $k = 1,\ldots,N$. Then we have $s^*_d \nearrow p^*$ as $d \rightarrow \infty$.
\end{proposition}

\begin{proof} Because any point $(\mu_1,\ldots,\mu_N)$ that is feasible for the optimization problem (\ref{equ:OptProbInvMeasure2}) induces a feasible point $(\mathbf{y}^{(1)},\ldots,\mathbf{y}^{(N)})$ for (\ref{equ:OptLinearFunctionalTruncated}) by its truncated moment sequence it follows $s^*_d \leq p_2^* = p^*$. Further, the feasible set is shrinking when $d$ increases, hence we have $s^*_{d+1} \geq s^*_d$. Convergence to $p^*$ follows by classical arguments as in \cite{MilanInvariantMeasure}. We will only sketch the proof, details can be found in \cite{MilanInvariantMeasure}. We start with a sequence $(\mathbf{y}^{(1,d)},\ldots,\mathbf{y}^{(N,d)})$ of optimal points for the truncated optimization problem (\ref{equ:OptLinearFunctionalTruncated}). As in \cite{MilanInvariantMeasure} we see that for each $k = 1,\ldots,N$ the components of the moment sequence $(\mathbf{y}^{(k,d)})_{d \in \N}$ are bounded. Extracting a component wise convergent subsequence leads to bounded sequences $\mathbf{y}^{(1)},\ldots,\mathbf{y}^{(N)}$. It remains to show that these sequences $\mathbf{y}^{(1)},\ldots,\mathbf{y}^{(N)}$ are moment sequences. Since taking the component wise limit preserves the property $M_d(g_i^{(k)}\mathbf{y}^{(k)}) \succeq 0$ for all $d \in \N$, $k = 1,\ldots,N$ and $i = 1,\ldots,l_k$ as well as $M_d(g_0\mathbf{y}^{(k)}) \succeq 0$ by Putinar's positivstellensatz this guarantees that each $\mathbf{y}^{(k)}$ is the moment sequence of a unique measure $\mu_k$ on $\Pi_{I_k}(X)$ for $k = 1,\ldots,N$. For the cost we write $G_k =  \sum\limits_{\alpha} G_{k,\alpha} x^\alpha$ -- where the sum is finite -- we have
\begin{equation}\label{equ:LimitCost}
	s^*_d = \sum\limits_{k = 1}^N  l_{\mathbf{y}^{(k)}}^d(G_k)= \sum\limits_{k = 1}^N  \sum\limits_{\alpha} G_{k,\alpha} \mathbf{y}^{(k,d)}_\alpha \rightarrow \sum\limits_{k = 1}^N  \sum\limits_{\alpha}G_{k,\alpha} \mathbf{y}^{(k)}_\alpha = \sum\limits_{k = 1}^N\langle G_k,\mu_k\rangle.
\end{equation}
Also the unit mass condition $\mathbf{y}{(k,d)}_0 = 1$, the invariance condition $\ell^d_{\mathbf{y}^{(k,d)}}(x^\alpha \circ f) = \ell^d_{\mathbf{y}^{(k,d)}}(x^\alpha)$ and the compatibility condition $\mathbf{y}^{(k,d)}_\alpha = \mathbf{y}^{(l,d)}_\alpha$ for $\alpha$ with $\alpha \in \N^{I_k \cap I_l}$ for $k,l = 1,\ldots,N$ pass to the limit and we get that the corresponding measures $\mu_1,\ldots,\mu_N$ are feasible for (\ref{equ:OptProbInvMeasure2}). Finally with (\ref{equ:LimitCost}) we get $p^* = p^*_2 \leq \langle G,\mu\rangle = \lim\limits_{d \rightarrow \infty} s^*_d$.
\end{proof}

\begin{remark}
	The optimization problem (\ref{equ:OptProbInvMeasure1}) has the advantage that its solution gives an optimal invariant measure while (\ref{equ:OptProbInvMeasure2}) respectively (\ref{equ:OptLinearFunctionalTruncated}) only provide (pseudo) moments of the marginals of such a measure. That is why, in case one is interested in the extremal measure and not only in the optimal value or its marginals, we suggest using (\ref{equ:OptProbInvMeasure1}) and to exploit its block structure. Another way is to solve (\ref{equ:OptProbInvMeasure2}) by means of solving (\ref{equ:OptLinearFunctionalTruncated}) and extend the obtained marginals to an (invariant) measure by solving a feasibility problem for linear matrix inequalities, i.e. formulate the problem: Find an invariant probability measure $\mu \in \mathcal{M}(X)_+$ with marginals $(\Pi_{I_k})_\# \mu$ given by the solutions of (\ref{equ:OptProbInvMeasure2}) respectively (\ref{equ:OptLinearFunctionalTruncated}). The invariance and marginal constraints are linear. This can be formulated in terms of moments and as in (\ref{equ:OptLinearFunctionalTruncated}) the non-negativity constraint can be expressed as linear matrix inequalities. Alternatively to enforce the invariance condition ergodic averages can be used after finding a probability measure with the given marginals, as mentioned in Remark \ref{rem:AveragingInvariantMeasure}. But those methods are still computationally complex because the complexity is inherited from the non-trivial task of finding measures with given marginals (and additional properties, as invariance for instance).
\end{remark}

An example where it is straight forward to give a measure $\mu$ with given marginals $\mu_k = (\Pi_{I_k})_{\#}\mu$ is when the marginals are atomic measures. A measure $\mu_k$ is atomic if it has a particle representations
\begin{equation}\label{equ:ParticleRep}
\mu_k = \sum\limits_{i = 1}^{m_k} a_{k,i} \delta_{x_{k,i}}
\end{equation}
for $m_k \in \N$, weights $a_{k,i} \in [0,\infty)$ and particle positions $x_{k,i} \in \R^{|I_k|}$. Invariant measures are typically not atomic, but several numerical methods provide approximations of measures by means of atomic measures. This, for instance, is the case for particle flow methods from \cite{ChizatBach}.
\, For this case we describe a natural choice for $\mu$, namely we will see that the particles $x_{k,i}$ used in the representations of $\mu$ are exactly the projections $\Pi_{I_k}(x_i)$ of some points $x_i$ in the whole state space. To do so we note the following:

\begin{enumerate}
	\item For two particle representations of a measure $\nu$, i.e. $\nu = \sum\limits_
{i = 1}^{m} a_{i} \delta_{x_{i}} = \sum\limits_{j = 1}^l b_j \delta_{y_j}$ for non-vanishing weights $a_i,b_j$ and pairwise different positions $x_i$ respectively $y_j$ it follows $m = l$, $x_i = y_{\sigma(i)}$ and $a_i = b_{\sigma(i)}$ for all $i = 1,\ldots,m$ for some permutation of indices $\sigma$. In other words the particle representation is unique.
	\item For a projection $\Pi_I$ onto a subset $I$ of the coordinates we have $(\Pi_I)_\# \nu = \sum\limits_
{i = 1}^{m} a_{i} \delta_{\Pi_I(x_i)}$
\end{enumerate}

Combining these two results gives the following Lemma.

\begin{lemma}\label{lem:AtomGlue}
    Let $I_1,\ldots,I_N$ be a covering of $\{1,\ldots,n\}$. For $k = 1,\ldots,N$ assume measures $\mu_k \in \mathcal{M}(\Pi_{I_k}(X))$ have particle representation as in (\ref{equ:ParticleRep}),
    \begin{equation}
    \mu_k = \sum\limits_{i = 1}^{m_k} a_{k,i} \delta_{x_{k,i}}
\end{equation}
    with $m_k \in \N$ and $x_{k,i} \in \Pi_{I_k}(X)$, and, for simplicity, assume that for $i = 1,\ldots,m_k$ the points $\Pi_{I_k \cap I_l}(x_{k,i})$ are pairwise distinct for all $k,l \in \{1,\ldots,n\}$. Assume further that the compatibility condition $(\Pi_{I_k \cap I_l})_\# \mu_k = (\Pi_{I_k \cap I_l})_\# \mu_l$ for all $k,l = 1,\ldots,N$ are satisfied. Then we have $m_k = m_l =: m$, $a_{k,i} = a_{l,i}$ for all $k,l \in \{1,\ldots,n\}$ (up to permutation of indices) and there exist $m$ particles $x_1,\ldots,x_m \in \R^n$ such that for all $k= 1,\ldots,N$ we have $x_{k,i} = \Pi_{I_k} (x_i)$ (up to permutation of indices for the $x_{k,i}$). In other words the measure $\mu:= \sum\limits_{i = 1}^m a_i \delta_{x_i}$ satisfies $(\Pi_{I_k})_\# \mu = \mu_k$.
\end{lemma}

For the optimization problem (\ref{equ:OptProbInvMeasure2}) that means: If (numerical) solutions $\mu_k$ of (\ref{equ:OptProbInvMeasure2}) have particle representations  (\ref{equ:ParticleRep}), we get, by Lemma \ref{lem:AtomGlue} from the above, that there exists a measure $\mu$ with particle representation and marginals $\mu_k$. Further, (numerical) invariance of $\mu$ is maintained, i.e. $\mu$ is feasible for (\ref{equ:OptProbInvMeasure2}), and $\mu$ produces the same cost, i.e. $\mu$ is optimal (by Proposition \ref{prop:EquivalentExtremalMeasures}).

We illustrate this by giving a numerical approximation of the global attractor. A typical approach is to just evolve a particle cloud for long time and then displaying the final section of the evolution of the point cloud. This can be interpreted as an approximation of an invariant measure by atomic measures.

The same procedure but based on subsystems can be interpreted in the following way: We choose, as above, an atomic representation/approximation of (invariant) measure $\mu$ with marginals given by the particle representations from the point clouds of the subsystems.

\begin{figure}[!h]
\begin{picture}(250,220)
\put(10,60){\includegraphics[width=78mm]{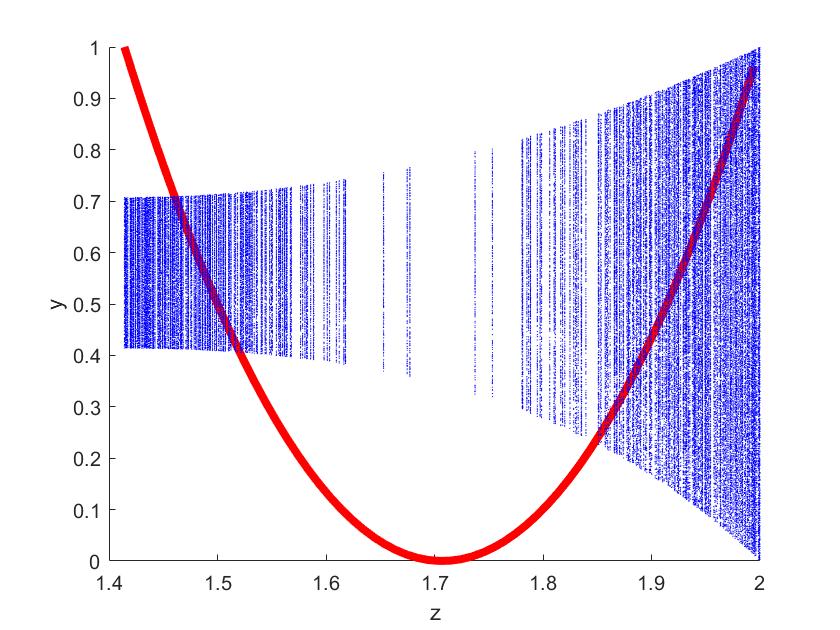}}
\put(220,60){\includegraphics[width=78mm]{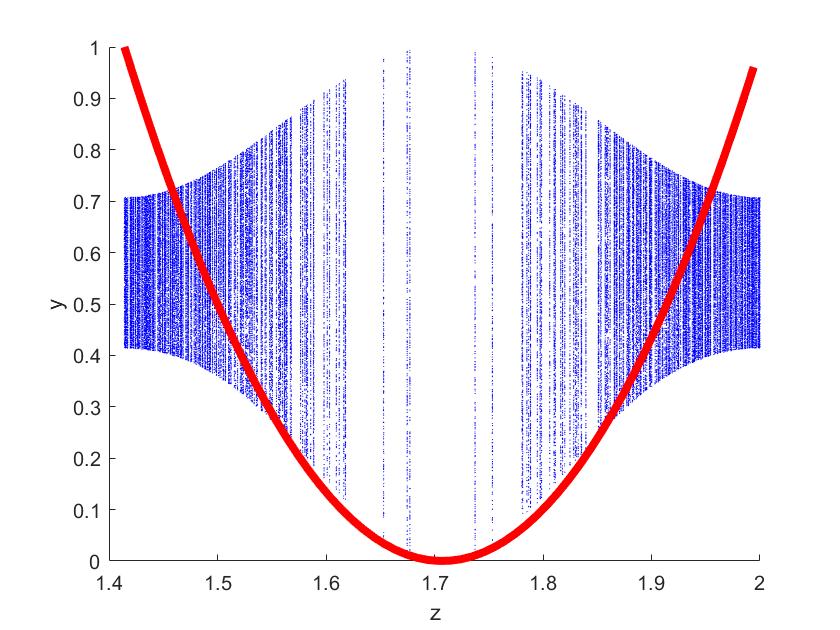}}
\end{picture}
\vspace{-15mm}
\caption{\footnotesize{Numerical particle representation of an invariant measure for (\ref{equ:TentMapCoupled}) for the subsystems where the samples of $z$ are drawn independently with respect to the (scaled) density displayed in red. The 550 $z$ samples coincide for both subsystems, the dynamics on the attractor were generated by evolving the initial value $x = 0.48934$ and $y = 0.8979573$ for 200 iterations and then plotting 300 evolution steps corresponding to each of $z$ samples.}}
\label{fig:Glue1}
\end{figure}

\begin{figure}[!h]
\begin{center}
\includegraphics[width=130mm]{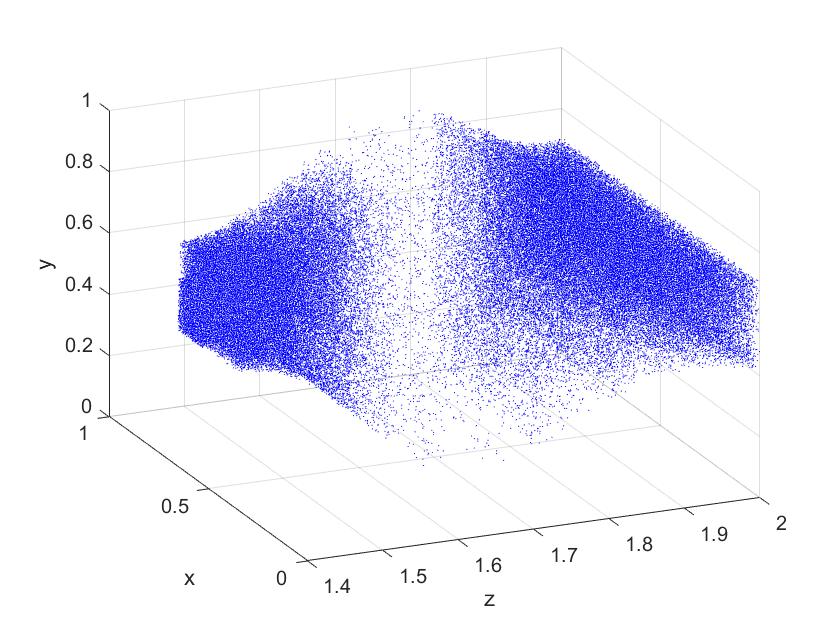}
\end{center}
\caption{\footnotesize{Displays the particle measures from Figures \ref{fig:Glue1} glued together along their common $z-$marginal, hence numerically approximating an invariant measure for the whole system (\ref{equ:TentMapCoupled}).}}
\label{fig:Glue2}
\end{figure}

We demonstrate this at the following artificial example in Figures \ref{fig:Glue1} and \ref{fig:Glue2} where we consider the coupled one dimensional tent maps $f_z$ for $z \in (1,2]$ defined by
\begin{equation}\label{equ:tentmap}
	f_z:[0,1] \rightarrow [0,1], x \mapsto \begin{cases} z x, & 0 \leq x \leq \frac{1}{2}\\
															z (1-x), & \frac{1}{2} \leq x \leq 1.
											  \end{cases}
\end{equation}
For fixed $z \in (\sqrt{2},2]$ the global attractor for the discrete dynamical system $x_{k+1} = f_z(x_k)$ is given by $[z(1-\frac{z}{2}),\frac{z}{2}]$. We consider the coupled system
\begin{eqnarray}\label{equ:TentMapCoupled}
	z_{k+1} & = & z_k\notag\\
	x_{k+1} & = &\begin{cases} \sin(\pi z_k) x_k, & 0 \leq x \leq \frac{1}{2}\\
															\sin(\pi z_k) (1-x_k), & \frac{1}{2} \leq x \leq 1.
											  \end{cases}\\
	y_{k+1} & = & \begin{cases} \left(\sqrt{2} + \left(\frac{z_k - \sqrt{2}}{2-\sqrt{2}}\right)(2-\sqrt{2})\right)  x_k, & 0 \leq x \leq \frac{1}{2}\\
															 \left(\sqrt{2} + \left(\frac{z_k - \sqrt{2}}{2-\sqrt{2}}\right)(2-\sqrt{2})\right) (1-x_k), & \frac{1}{2} \leq x \leq 1.
											  \end{cases}\notag\\\notag
\end{eqnarray}
For each state $x$ and $y$ we transformed the parameter $z_k$ in order to introduce different behaviour in their dynamics. The dynamics (\ref{equ:TentMapCoupled}) is of the form (\ref{equ:Sparsex0x1x2}). Hence it decouples into the corresponding subsystems induced by $I_1 = \{1,2\}$ and $I_2= \{1,3\}$, by denoting $x^1 = z$, $x^2 = x$ and $x^3 = y$.

The Figures \ref{fig:Glue1} and \ref{fig:Glue2} are generated by drawing 550 random samples of $z$ according to a distribution with (scalled) density displayed in red in Figures \ref{fig:Glue1} and \ref{fig:Glue2}. For each such sample $z$ we evolve an initial state $x$ respectively $y$ according to (\ref{equ:TentMapCoupled}) for 500 steps. Figure \ref{fig:Glue1} displays in blue the last 300 evolutions of the sates $x$ respectively $y$ for each sample $z$ -- this gives an atomic approximation of an invariant measure for the subsystems. The two measures represented in Figure \ref{fig:Glue1} are ``glued" together, as described above, in Figure \ref{fig:Glue2}, by using the particles from the two subsystems connected by their common state samples $z$. Note that this gives nothing else than exactly the point cloud given by evolving the initial states $x,y$ on the samples $z$ simultaneously for the whole system.


\begin{remark}
	A decoupling into subsystems allows to treat the subsystems independently which allows to take specific properties of the subsystems into account. In the case of EDMD, different observables can be used for the different subsystems based on different properties of the subsystems. That includes more observables for one subsystem, for example if this subsystem shows more complexity, or a different class of observables as for example trigonometric polynomials, piecewise linear functions etc. Similarly for the computation of invariant measures; for the subsystems different observables/test functions can be chosen or if the invariant measure is computed using particles, i.e. point masses, different numbers of particles can be used for different subsystems. 
\end{remark}

\section{Conclusion}
We have presented a framework for exploiting specific sparse structures, namely subsystems, of dynamical systems for the Koopman operator and Perron-Frobenius operators. The notion of subsystems, that we use in this text, allows to decompose the dynamical system and, as a consequence, the Koopman operator decomposes into Koopman operators on the subsystems. 
In contrast to the notion of subsystems in \cite{ElkinSubsystem}, the subsystems we treat describe dynamic independence of several states rather than functional dependence of different states as in \cite{ElkinSubsystem}. That makes available simple algorithms from graph theory, based on the sparsity graph of the dynamics, in order to find the subsystems we are interested in. But at the same time we neglect functional correlation between dependent states. Exploiting term-sparsity structures for dynamical systems can be found in \cite{JieDS} and considering such structures for the Koopman operator could provide a future perspective complementing our work.

With respect to data science, the advantage of our approach compared to other sparsification techniques based on data is that we make use of the inherited (exact) sparse structure of the dynamical system. Hence for such cases the sparse structure in the data comes naturally and allows lower dimensional treatments of the system. The presented results in this work show that objects such as (principal) eigenfunctions and invariant measures can be recovered from lower dimensional subsystems and hence allow a reduction in complexity.

This approach is related to sparsity patterns in the choice of observables, or limited information to put it in different words. An approach that incorporates this can be found in \cite{Aqip}. In contrast to \cite{Aqip} we treat the case where we know a-priori that restricting to specific observables (i.e. projections onto the subsystems) allows a full recovery of the system.

Our approach can also be applied to system identification via the Koopman operator (\cite{SystemIdentificationMauroy}, \cite{RosenfeldSystemIdentification}) and combined with the symmetry approach from \cite{Symmetry}.

The computational gain of our approach depends highly on the factorization into subsystems, hence a good choice for these can be highly beneficial and is treated in \cite{SparsePaper}.


\section{Acknowledgements}
This work has been supported by European Union’s Horizon 2020 research and innovation programme under the Marie Skłodowska-Curie Actions, grant agreement 813211 (POEMA), by the Czech Science Foundation (GACR) under contract No. 20-11626Y and by the AI Interdisciplinary Institute ANITI funding, through the French “Investing for the Future PIA3” program under the Grant agreement n$^\circ$ ANR-19-PI3A-0004.

The first author wants to thank Matthew D. Kvalheim for an inspiring discussion about principal eigenfunctions which has strongly influenced the corresponding part in this text.

\bibliographystyle{siamplain}
\bibliography{references}

\end{document}